\documentclass[12pt]{amsart}
\usepackage[dvips]{graphicx}
\newtheorem{thm}{Theorem}[section]

\newtheorem{lem}[thm]{Lemma}
\newtheorem{prop}[thm]{Proposition}
\newtheorem{defn}[thm]{Definition}
\newtheorem{note}[thm]{Note}

\newtheorem{example}[thm]{Example}
\newtheorem{lemma}[thm]{Lemma}
\newtheorem{proposition}[thm]{Proposition}
\newtheorem{corollary}[thm]{Corollary}
\newtheorem{definition}[thm]{Definition}
\newcommand{\Eop}{\ \textrm{q.e.d.}}
\newcommand{\Proof}[1]{\noindent {\itshape  #1}}
\newcommand{\iIR}{\mathbb{R}}
\newcommand{\iIZ}{\mathbb{Z}}
\newcommand{\prd}[1]{(\!\!\!(#1)\!\!\!)}

\begin{document}

\title[K$(\pi, 1)$'s for Artin Groups]
        {K$(\pi, 1)$'s for Artin Groups of Finite Type}

\author[T.~Brady]{Thomas Brady}

        \address{School of Mathematical Sciences\\
                Dublin City University\\
                Glasnevin, Dublin 9\\
                Ireland}
        \email{tom.brady@dcu.ie}
\author[C.~Watt]{Colum Watt}

        \address{School of Mathematics\\
                Trinity College\\
                Dublin 2\\
                Ireland}
        \email{colum@maths.tcd.ie}



\maketitle

\section{Introduction.}
This paper is a continuation of a programme to construct new K$(\pi,1)$'s for Artin groups of 
finite type which began in \cite{brady} with Artin groups on $2$ and $3$ generators and 
was extended to braid groups in \cite{brady1}.  These K$(\pi,1)$'s differ from those in 
\cite{davis} in that their universal covers are simplicial complexes.
\vskip .2cm
In \cite{brady} a complex is constructed 
whose top-dimensional cells correspond to minimal factorizations of a Coxeter element as 
a product of reflections in a finite Coxeter group.  Asphericity is established in low 
dimensions using a metric of non-positive curvature.   Since the non-positive curvature 
condition is difficult to check in higher dimensions a combinatorial approach is used in 
\cite{brady1} in the case of the braid groups.  
\vskip .2cm
It is clear from \cite{brady1} that the 
techniques used can be applied to any finite Coxeter group $W$.  
When $W$ is equipped with the partial order given by reflection length and 
$\gamma$ is a Coxeter element in $W$, the construction of the K$(\pi,1)$'s is exactly 
analogous provided that the interval $[I, \gamma]$ forms a lattice.  In dimension $3$, see 
\cite{brady}, establishing this condition amounts to observing that two planes through the 
origin meet in a unique line.  In the braid group case, see \cite{brady1},  
where the reflections are transpositions and the Coxeter element is an $n$-cycle this 
lattice property is established by identifying $[I, \gamma]$ with the lattice of 
non-crossing partitions of $\{1,2,\dots ,n\}$.
\vskip .2cm
In this paper, we consider the Artin groups of type $C_n$ and $D_n$.  Thus, for each 
finite reflection group $W$ of type $C_n$ or $D_n$, partially ordered by reflection length, we 
identify a lattice inside $W$ and use it to construct a finite aspherical complex $K(W)$.
  In the $C_n$ case this lattice coincides with the lattice of noncrossing 
partitions of $\{1,2,\dots , n, -1, \dots, -n\}$ studied in \cite{reiner}.  
The final ingredient is to prove that $\pi_1(K(W))$ is isomorphic to $A(W)$, the associated 
finite type Artin group.  As in \cite{brady} and \cite{brady1} this involves a lengthy 
check that the obvious maps between the two presentations are well-defined. 
\vskip .2cm
David Bessis has independently obtained similar results which can be seen at \cite{bessis}.  
His approach exploits in a clever way the extra structure given by viewing these 
groups as complex reflection groups.  In addition, he has verified that in the 
exceptional cases that the interval $[I, \gamma]$ forms a lattice and that the 
corresponding poset groups are isomorphic to the respective Artin groups of finite type.  
Combined with the results of our section $5$ below this provides the new  K$(\pi,1)$'s 
in these cases and we thank him for drawing our attention to this fact.

\vskip .2cm
In section 2 we collect some general facts about the reflection length function on finite 
reflection groups and the induced partial order.   In section 3 we study the cube group 
$C_n$ and its index two subgroup $D_n$. 
In section 4 we identify the subposets of interest in $C_n$ and $D_n$ and show that 
they are lattices.  
In section 5 we define the poset group $\Gamma(W, \alpha)$ associated to the 
interval $[I, \alpha]$ for $\alpha \in W$.
In the case where $[I, \alpha]$ is a lattice we construct the complexes $K(W, \alpha)$ and 
show that they are K$(\pi, 1)$'s.   Section 6 shows that the groups $\Gamma(C_n, \gamma_C)$ and 
$\Gamma(D_n, \gamma_D)$ are indeed the Artin groups of the 
appropriate type when $\gamma_C$ and $\gamma_D$ are the respective Coxeter elements.

\section{A partial order on finite reflection groups.}

Let $W$ be a finite reflection group with reflection set ${\mathcal R}$ and identity element 
$I$.    We let $d:W \times W \to {\bf Z}$ be the distance function in the Cayley graph of 
$W$ with generating set ${\mathcal R}$ and   
define the \emph{reflection length function}  
$l: W \to {\bf Z}$ by $l(w) = d(I, w)$.  
So $l(w)$ is the length of the shortest product of reflections yielding the element $w$. 
It follows from the triangle inequality for $d$ that 
$ l(w) \le l(u) + l(u^{-1}w) $
for any $u,w \in W$.  
\vskip .2cm
\begin{defn}\label{d:order}
We introduce the relation $\le$ on $W$ by declaring 
\[u \le w \quad \quad \Leftrightarrow \quad \quad l(w) = l(u) + l(u^{-1}w).\]
\end{defn}
\vskip .2cm  
Thus $u \le w$ if and only if  there 
is a geodesic in the Cayley graph from $I$ to $w$ which 
passes through $u$.  Alternatively, equality occurs if and only if there is a shortest 
factorisation of $u$ as a product of reflections which is a prefix of a shortest 
factorisation of $w$.
It is readily shown that $\le$ is reflexive, antisymmetric and 
transitive so that $(W, \le)$ becomes a partially ordered set.   
\vskip .2cm
Since $(u^{-1}w)^{-1}w = w^{-1}uw $ is conjugate to $u$ it follows that $u^{-1}w \le w$ 
whenever $ u \le w$.  Furthermore, whenever $\alpha \le \beta \le \gamma$ we have 
\[l(\gamma) = l(\alpha) + (l(\alpha^{-1}\beta)+l(\beta^{-1}\gamma)),\]
so that $\alpha^{-1}\beta \le \alpha^{-1}\gamma$.
\vskip .2cm
We recall some general facts about orthogonal transformations from \cite{bradywatt}.  If 
$A \in O(n)$, we associate to $A$ two subspaces of ${\bf R}^n$ , namely
\[M(A) = \mbox{im}(A-I) \quad \mbox{ and } \quad F(A) = \mbox{ker}(A-I).\]
We recall that $M(A)^{\perp} = F(A)$.  We use the notation $|V|$ for $\dim(V)$ when $V$ is a 
subspace of ${\bf R}^n$.   It is shown in \cite{bradywatt} that 
\[|M(AC)| \le |M(A)| +|M(C)|\]
 We define a partial order on $O(n)$ by 
\[A \le_o B \quad \Leftrightarrow \quad |M(B)| = |M(A)| + |M(A^{-1}B)|\]
and we note that $A \le_o B$ if and only if $M(B) = M(A) \oplus M(A^{-1}B)$.   In particular 
$A \le_o B$ implies that $M(A) \subseteq M(B)$ or equivalently $F(B) \subseteq F(A)$.  
The main result we will use from \cite{bradywatt} is that for each $A \in O(n)$ and each 
subspace $V$ of $M(A)$ there exists a unique $B \in O(n)$ with $B \le_o A$ and $M(B) = V$.  
\vskip .2cm
Our finite reflection group $W$ is a subgroup of $O(n)$, so the 
results of \cite{bradywatt} can be applied to the elements of $W$.  We begin with a 
geometric interpretation of the length function $l$ on $W$.  
\begin{proposition}  $l(\alpha) = |M(\alpha)| = n-|F(\alpha)|$, for $\alpha \in W$.
\end{proposition}
\Proof{Proof.} First note that the proposition holds when $\alpha = I$ so we will 
assume $\alpha \ne I$ and let $ k = |M(\alpha)| > 0$.
\vskip .2cm
To establish the inequality $l(\alpha) \le k$ we  show that $\alpha$ can be 
expressed as a product of $k$ reflections.  We will use induction 
on $k$ noting that the case $k = 1$ is immediate.  Consider the 
subspace $F(\alpha) \ne \iIR^n$.  Recall from part (d) of Theorem 1.12 of 
\cite{humphreys} that the subgroup $W'$ of $W$ of elements which fix $F(\alpha)$ 
pointwise is generated by those reflections $R$ in $W$ satisfying $F(\alpha) \subset F(R)$.  
Since $\alpha \ne I$ there exists at least one such reflection $R$.  Since 
$M(A) = F(A)^{\perp}$ we have $M(R) \subset M(\alpha)$.  The unique orthogonal transformation 
induced on $M(R)$ by $\alpha$ must be $R$ by Corollary~3 of \cite{bradywatt}.  Hence 
$R \le_o \alpha$ and 
\[|M(R\alpha)| = |M(\alpha)|- |M(R)| = k - 1.\]
By induction $R\alpha$ can be expressed as a product of $k - 1$ reflections and hence there 
is an expression $\alpha = R_1\dots R_k$ for $\alpha$ as a product of $k$ reflections.  
We note that by construction each of these reflections $R_i$ satisfies 
$M(R_i) \subset M(\alpha)$.
\vskip .2cm
To establish the other inequality suppose $\alpha = S_1S_2\dots S_m$ is an expression 
for $\alpha$ as a product of $m$ reflections realizing $l(\alpha) = m$.  Repeated use of 
the identity $|M(AC)| \le |M(A)| +|M(C)|$ gives 
\[k = |M(\alpha)| \le |M(S_1)| + \dots + |M(S_m)| = m = l(\alpha). \quad \quad \quad  \Eop\]
\vskip .2cm
In particular the partial order $\le$ on $W$ is a restriction 
of the partial order $\le_o$ on $O(n)$ and we will drop the subscript from $\le_o$ from 
now on.  The following lemma is immediate.
\begin{lemma}
\label{sublength}
Let $W$ be a finite Coxeter group with reflection set ${\mathcal R}$
and let $W_1$ be a subgroup generated by a subset ${\mathcal R}_1$ of ${\mathcal R}$.
Then the length function for $W_1$ is equal to the restriction
to $W_1$ of the length function for $W$.
\end{lemma}
\begin{definition}
For each  $\delta \in W$ 
we define the reflection set of $\delta $, $S_{\delta}$, by
$S_{\delta} = \{R \in {\mathcal R} \mid r \le \delta\}$.
\end{definition}
Repeated application of 
$A \le B \Rightarrow  |M(B)| = |M(A)| + |M(A^{-1}B)|$ gives 
$M(\delta ) = \mbox{Span}\{M(R)\mid R \le \delta\}$ so that $S_\delta$ determines $M(\delta)$.  
However, in the case where $\delta \le \gamma$, $\delta$ itself is determined by $\gamma$ and 
$S_{\delta}$ since $\delta$ is the unique orthogonal transformation induced on $M(\delta)$ 
by $\gamma$.  The following results are consequences of this fact.
\begin{lemma}\label{orderchar}
If $\alpha, \beta \le \gamma$ in $W$ and $S_{\alpha} \subseteq S_{\beta}$ then 
$\alpha \le \beta$.
\end{lemma}
\Proof{Proof.}  $M(\alpha) \subset M(\beta) \subset M(\gamma)$ and by uniqueness 
the transformation induced on $M(\alpha)$ by $\beta$ is the same as the transformation 
induced by $\gamma$, namely $\alpha$. \hfill \Eop
\begin{lemma}\label{meetchar}
Suppose  $\alpha, \beta \le \gamma$ in $W$. 
If there is an element $\delta \in W$ with $\delta \le \gamma$ and 
$ S_{\delta} = S_{\alpha} \cap S_{\beta}$ 
then  $\delta $ is the greatest lower bound of $\alpha$ and $\beta$ in $W$, that is, 
if $\tau \in W$ satisfies $\tau \le \alpha, \beta$ then $\tau \le \delta$.
\end{lemma}
\section{The Cube groups $C_n$ and $D_n$.}
For general facts about the groups $C_n$ and $D_n$ see \cite{bourbaki} or \cite{humphreys}.  
Let $I = [-1,1]$ and let 
$C_n$ denote the group of isometries of the cube $I^n$ in $\iIR^n$.
That is
\[ C_n = \{ \alpha \in O(n) : \alpha(I^n) = I^n \} \]
Let $e_1,\ldots,e_n$ denote the standard basis for $\iIR^n$
and 
let $x_1,\ldots,x_n$ denote the corresponding coordinates.
The set ${\mathcal R}_c$ of all reflections in $C_n$ consists of
the following $n^2$ elements. 
For each $i = 1,\ldots,n$, reflection in the hyperplane $x_i=0$ 
is denoted $[i]$ and also by $[-i]$.
For each $i\neq j$, reflection in the hyperplane $x_i=x_j$
                is denoted by any one of the four expressions
                $\prd{  i,j  } $, $\prd{  j,i  } $,
                $\prd{  -i,-j } $ and $\prd{  -j,-i } $,
                while reflection in the plane $x_i= -x_j$
                is denoted by any one of the four
                expressions $\prd{  i,-j  } $,
                 $\prd{  -i,j  } $,
                 $\prd{  j,-i  } $, and 
                 $\prd{  -j,i  } $.
The set of these $n(n-1)$ reflections, in hyperplanes of the form $x_i = \pm x_j$, 
is denoted ${\mathcal R}_d$ and the subgroup they generate, $D_n$, is well known 
to be an index two subgroup of $C_n$.
The group $C_n$ acts 
on the set $\{e_1,\ldots,e_n,-e_1,\ldots, -e_n\}$
in the obvious manner and this action satisfies
$ \alpha\cdot (-e_i) = -(\alpha\cdot  e_i)$   
for each $i$ and each $\alpha \in C_n$.
Thus we obtain an injective homomorphism
$ p$ from   $ C_n$ into the group $\Sigma_{2n}$
of permutations  of the set
$ \{ 1,2,\ldots,n,-1,-2,\ldots,-n\}$.
Note that for each  $i$, $p([i])$
is a transposition in $\Sigma_{2n}$,
while each  element of ${\mathcal R}_d$
is mapped to a product of two disjoint transpositions.
Thus $p(D_n)$ is contained in the subgroup of even 
permutations.

For each cycle  $c = (i_1,\ldots,i_r)$ in $\Sigma_{2n}$,
we define the cycle $\bar{c}$ by
\[\bar{c} = (-i_1,\ldots,-i_r) \]
Note that $\bar{c} = z_0cz_0$ where $z_0 = (1,-1)(2,-2)\ldots(n,-n)$
has order two. Note also that $z_0= p(\zeta_0)$ where
$\zeta_0 = [1][2]\cdots[n]$ is the nontrivial element in the 
centre of $C_n$.

\begin{proposition}
\label{image}
The image~$p(C_n)$ is the 
centraliser~$Z(z_0)$ of $z_0$ in $\Sigma_{2n}$.
It consists of all products of disjoint cycles of the form
\begin{equation}
\label{form}
 c_1\bar{c}_1\ldots c_k\bar{c}_k\gamma_1\ldots\gamma_r 
\mbox{ \ where \ } \gamma_j = \bar{\gamma}_j \ \ \ \forall \  j=1,\ldots,r. 
 \end{equation}
The image $p(D_n)$ consists of all elements of the form~(\ref{form})
with $r$ even.
\end{proposition}

\Proof{Proof.}
Since $z_0$ has order $2$ and
$z_0 c_1c_2\ldots c_k z_0 = \bar{c}_1\bar{c}_2\ldots\bar{c}_k$
for any product of cycles in $\Sigma_{2n}$, 
it follows that
the centraliser~$Z(z_0)$  consists of those products
of disjoint cycles
$c_1c_2\ldots c_k$ for which
\[ c_1c_2\ldots c_k = \bar{c}_1\bar{c}_2\ldots\bar{c}_k\]
By uniqueness (up to reordering) of cycle decomposition
in $\Sigma_{2n}$, for each $i$ either $c_i = \bar{c}_j$
for some $j \neq i$ or else $c_i=\bar{c}_i$.
It follows that the centraliser of $ z_0 $ is
precisely the set of elements in $\Sigma_{2n}$ of the form~(\ref{form}).

For each $\alpha \in C_n$, the identity $\zeta_0 \alpha \zeta_0 = \alpha$
implies that $p(\alpha)$ lies in the centraliser of $z_0$.
Thus $p(C_n)\subset Z(z_0)$.
In the reverse direction,
if $c= (i_1,\ldots,i_k)$ is disjoint from $\bar{c}$,
one may readily verify that
\begin{equation}
\label{pfact}
c\bar{c}
= p\left( \prd{i_1,i_{2}} \, \prd{i_2,i_{3}}  
                     \ldots \prd{  i_{q-1},i_{q }}  \right) 
\end{equation}
Likewise, if $c = \bar{c}$ then $c$ must be 
the form $c = (i_1,\ldots,i_k,-i_1,\ldots,-i_k)$ for 
some $-n \leq i_1,i_2, \ldots,i_k \leq n$ 
and one may verify that
\begin{eqnarray}
\label{bfact}
c &=& (i_1,-i_1)(i_1,i_2)(-i_1,-i_2)\ldots(i_{k-1},i_k)(-i_{k-1},-i_k) \\
 &=& p\left( [i_1]\prd{i_1,i_2}  
                \ldots\prd{i_{k-1},i_k}  \right) 
\end{eqnarray}
It follows that 
any element of the form~(\ref{form}) lies in $p(C_n)$
and hence $p(C_n) =  Z(z_0)$.

Let $\alpha \in D_n$ and write 
$p(\alpha)=c_1\bar{c}_1\cdots c_a\bar{c}_a
\gamma_1\cdots \gamma_b$.
Since $p(\alpha)$ and each $c_i\bar{c}_i$ is an even permutation 
while each $\gamma_j$ is an odd permutation, $r$ must be even.
To show that every element of the form~(\ref{form}) 
with $r$ even is in $p(D_n)$, we need only
note the following facts.
\begin{itemize}
\item If the cycle~$c$ is disjoint from $\bar{c}$ then equation~(\ref{pfact})
implies that $c\bar{c} \in p(D_n)$.
\item If $i \neq j$ then $[i][j] = \prd{i,j}\,\prd{i,-j}$ and hence is
an element of $p(D_n)$. It now follows from equation~(\ref{bfact})
that if $c_1=\bar{c}_1$ and $c_2=\bar{c}_2$ are disjoint cycles then
$c_1c_2 \in p(D_n)$.
\hfill $\Eop$
\end{itemize}

\bigskip

\noindent
\paragraph{\bf Notation}
From now on we will identify $C_n$ and $D_n$ with their
respective images in $\Sigma_{2n}$.
If a cycle $c = (i_1,\ldots,i_k)$ is disjoint from $\bar{c}$ then we 
write
\[ \prd{  i_1,\ldots,i_k  }  = c \bar{c} =
(i_1,\ldots,i_k)(-i_1,\ldots,-i_k) \]
and we call $c\bar{c}$ a \emph{paired cycle}. If $k=1$ then
$c = (i_1)$ and the paired cycle~$c\bar{c} = \prd{i_1}$ 
fixes the vector~$e_{i_1}$.
If $c = \bar{c} = (i_1,\ldots,i_r,-i_1,\ldots,-i_r)$
then we say that $c$ is a balanced cycle and we write
\[ c =  [i_1,\ldots,i_k] .\]
This notation is consistent with that 
introduced earlier for the elements of the generating set ${\mathcal R}_c$.
With these conventions, proposition~\ref{image} states 
that each element of $C_n$
may be written as a product of disjoint paired cycles and balanced cycles.
If $\alpha \in C_n$ fixes the standard basis vector~$e_i$
then we will assume that the paired cycle $\prd{i}$
appears in the corresponding expression~(\ref{form}) for $\alpha$.
\vskip .2cm
\noindent
Denote the length function for $C_n$ 
with respect to the generating set ${\mathcal R}_c$ 
by $l$. Lemma~\ref{sublength} allows us to use the same
symbol~$l$ for 
the length function of $D_n$ with respect to the set ${\mathcal R}_d$.
The length function for $\Sigma_{2n}$
with respect to the set~$T$ of all transpositions is
denoted by $L$.

\begin{lemma}
\label{fixedspaces}
The fixed space $F(\prd{i_1,\ldots,i_k})$
has dimension $n-k+1$ and is given by 
\[ \{x \in \iIR^n : x_{i_1}=x_{i_2}= \cdots = x_{i_k} \} \]
where $x_i$ means $-x_{|i|}$ for $i < 0$.
The fixed space $F([i_1,\ldots,i_k])$
has dimension $n-k$ and is given by
\[\{x \in \iIR^n : x_{i_1}=x_{i_2}= \cdots = x_{i_k} =0 \} \]
\end{lemma}

\Proof{Proof.} By inspection. \hfill $\Eop$

\begin{lemma}
\label{pairedlength}
The $l$-length of a paired cycle~$c\bar{c} = \prd{  i_1,\ldots, i_k  } $
is $k-1$. Moreover, no 
minimal length factorisation of $c\bar{c}$ as
a product of elements of ${\mathcal R}_c$ contains a generator
of the form $[i]$.
\end{lemma}

\Proof{Proof.}
The fixed space $F(c\bar{c})$ has dimension $n-k+1$ by lemma \ref{fixedspaces}
and thus $l(c\bar{c}) = n - (n-k+1) = k-1$. 

If a minimal $l$-length  factorisation of $c\bar{c}$
contained a term of the form $[i]$, we would obtain a factorisation of
$c\bar{c}$ as a product of fewer than $2(k-2) + 1 = 2k-3$
transpositions. As $L(c\bar{c}) = 2k-2$
this is impossible. \hfill $\Eop$
\begin{lemma}
\label{balancedlength}
The $l$-length of $\gamma  = [j_1,\ldots ,j_r ]$ as
a product of elements of ${\mathcal R}_c$ is $r$. Moreover
any minimal length factorisation of 
$\gamma$ as
a product of elements of ${\mathcal R}_c$ contains exactly one
generator of the form $[i]$.
\end{lemma}

\Proof{Proof.}
As the fixed space $F(\gamma)$ is $(n-r)$-dimensional
by lemma~\ref{fixedspaces}, we find $l(\gamma) = n-(n-r) = r$. 

As $L(\gamma) = 2r-1$, any factorisation of $\gamma$ as a
product of $r$ elements of ${\mathcal R}_c$ 
can contain at most one generator of the form $[i]$. 
If such a factorisation contained no element of this form, we would
have an expression for $\gamma$ as a product of an 
even number of transpositions. But this contradicts the fact
that the $2r$-cycle $\gamma$ has odd parity in $\Sigma_{2n}$.
\hfill $\Eop$

\begin{proposition}
\label{dimfix}
If $\alpha = c_1\bar{c}_1\ldots c_a\bar{c}_a \gamma_1\ldots\gamma_b \in C_n$
is a product of disjoint cycles
then 
$$l(\alpha) =\sum_{i=1}^a l(c_i\bar{c}_i) + \sum_{j=1}^b l(\gamma_j)$$
\end{proposition}

\Proof{Proof.}
By choosing a new basis from $\{e_1,\ldots,e_n,-e_1,\ldots,-e_n\}$
if necessary,
we may assume that $c_i = (j_{i-1}+1,j_{i-1}+2,\ldots,j_i)$
and $\gamma_i=[k_{i-1}+1,k_{i-1}+2,\ldots,k_i]$ where
$1=j_0<j_1<\cdots <j_a <j_a+1=k_0<k_1< \cdots < k_b=n$.
Then $c_i\bar{c}_i$ (resp.~$\gamma_j$) maps
$U_i = {\rm span}(e_{j_{i-1}+1},e_{j_{i-1}+2},\ldots,e_{j_i})$
(resp.~$V_i = {\rm span}(e_{k_{i-1}+1},e_{k_{i-1}+2},\ldots,e_{k_i})$)
to itself and leaves all the other $U$'s and $V$'s
pointwise fixed. As $c_i\bar{c}_i$ (resp.~$\gamma_j$)
fixes a $1$ (resp.~$0$) dimensional subspace
of $U_i$ (resp.~$V_j$), we see that $\alpha$ fixes
an $a$-dimensional subspace of $\iIR^n$.
Therefore $l(\alpha) = n-a$. Since
$\sum (1+l(c_i\bar{c}_i)) +  \sum l(\gamma_j) = n$ by
lemmas~\ref{pairedlength} and \ref{balancedlength}, the result follows. 
\hfill $\Eop$

\bigskip

Consider now the effect of multiplying $\alpha \in C_n$ 
on the right by a reflection~$R= \prd{i,j}$ or $R=[i]$.
It is clear that  only those cycles which contain an
integer of $R$ will be affected. The following example lists
the possibilities and the corresponding changes in lengths.

\begin{example}
\label{ex:comp}

The following four identities can be verified directly.
\begin{eqnarray*}
 [i_1, i_2, \ldots, i_k]{[i_k]} &=& {\prd{i_1,  i_2, \ldots , i_k}} \\
 {[i_1, i_2, \ldots, i_k]} \prd{i_j,i_k}  
&=& [ i_1, \dots, i_j] \prd{i_{j+1},  i_{j+2}, \ldots ,  i_k }  \\
 \prd{  i_1, i_2, \ldots, i_k }\,\prd{i_j,i_k}  &=&
 \prd{  i_1, \ldots, i_j }\, \prd{  i_{j+1},  i_{j+2}, \ldots , i_k } \\
 {[i_1,\ldots,i_j]}{[i_{j+1},\ldots,i_k]}\prd{ -i_j,i_k}  &=&
 \prd{i_1,i_2,\ldots,i_k } 
 \end{eqnarray*}
Since each reflection has order $2$, the following identities
are immediate.
\begin{eqnarray*}
{[i_1, i_2, \ldots , i_k]} &=& \prd{  i_1,  i_2, \ldots , i_k }[i_k] \\
{[i_1, i_2, \ldots , i_k]}&=&  {[ i_1, \ldots, i_j]} 
           \prd{i_{j+1},  i_{j+2}, \ldots ,  i_k }\,\prd{i_j,i_k}  \\
\prd{  i_1, i_2, \ldots, i_k }&=& \prd{  i_1, \ldots, i_j } 
             \,\prd{  i_{j+1},  i_{j+2}, \ldots , i_k }\,\prd{i_j,i_k}  \\
{[i_1,\ldots,i_j]}{[i_{j+1},\ldots,i_k]} &=& 
                 \prd{i_1,i_2,\ldots,i_k }   \,\prd{ -i_j,i_k}  
\end{eqnarray*} 
By proposition~\ref{dimfix}, we see that 
\begin{eqnarray*}
 l([i_1, i_2, \ldots, i_n]) &=& l(\prd{  i_1,  i_2, \ldots , i_n }) +1\\
 l([i_1, i_2, \ldots, i_n]) &=& l([ i_1, \ldots, i_j] 
                        \prd{i_{j+1},  i_{j+2}, \ldots ,  i_n }) +1 \\
 l(\prd{  i_1, i_2, \ldots, i_n }) &=& l(\prd{  i_1, \ldots, i_j } 
\,\prd{  i_{j+1},  i_{j+2}, \ldots , i_n })+1 \\
 l([i_1,\ldots,i_j][i_{j+1},\ldots,i_k])  &=& 
      l(\prd{i_1,i_2,\ldots,i_k } ) +1
\end{eqnarray*}
\end{example}
\begin{definition}
Let $\sigma = c_1c_2\cdots c_k$ and $\tau = d_1 d_2 \cdots d_l$
be two products of disjoint cycles in $\Sigma_{2n}$. We say that $\sigma$
is contained in $\tau$ (and write $\sigma \subset \tau$) if for each $i$
we can find $j$  such that the set of integers in the cycle~$c_i$
is a subset of the set of integers in the cycle $d_j$. This notion
restricts to give a notion of containment for elements of $C_n$.

A reflection $\prd{i,j}$ is s-contained in
$\alpha = c_1\bar{c}_1\ldots c_a\bar{c}_a \gamma_1\ldots\gamma_b \in C_n$
(and we write $\prd{i,j} \sqsubset \alpha$)
if $i$ is contained in $\gamma_k$ and $j$ is contained in $\gamma_l$
for some $k \neq l$.
\end{definition}
\begin{lemma}
\label{contained}
Let $\alpha \in C_n$ and $R \in \mathcal{R}_c$.
Then $R \le \alpha$ if and only if
$R \subset \alpha$ or $R \sqsubset \alpha$.
\end{lemma}

\Proof{Proof.} By proposition~\ref{dimfix} and the calculations in example~\ref{ex:comp}
we see that $l(\alpha R) < l(\alpha)$ if and only if $R\subset \alpha$
or $R \sqsubset \alpha$. Since $R \leq \alpha$ if and only if
$l(\alpha R) < l(\alpha)$, the lemma follows. \hfill $\Eop$

\section{The lattice property}
In this section we show that the interval $[1, \gamma]$ in $(W \le)$ is a lattice for 
$W = C_n, D_n$ and $\gamma$ a Coxeter element in $W$.  
Since all Coxeter elements in $W$ are 
conjugate we can choose our favourite one in each case.
\begin{definition}  We choose the Coxeter  elements $\gamma_C$ in $C_n$ and 
$\gamma_D$ in $D_n$ given by $\gamma_C = [1,2,\dots , n]$ and $\gamma_D = [1][2,3,\dots ,n]$.
\end{definition}

\begin{proposition}
\label{b_i}
Write the Coxeter element~$\gamma_C \in C_n$ (resp.~$\gamma_D \in D_n$)
as 
$\gamma_C = R_1R_2\ldots R_n$ (resp.~$\gamma_D = R_1R_2\ldots R_n$)
for  reflections $R_1,\ldots,R_n$ in $\mathcal{R}_c$
(resp.~$\mathcal{R}_d$) and let $b_i$ denote the number
of balanced cycles in $R_1R_2\cdots R_i$.
Then there exists $i_0$ such that
$b_i = 0$ for $i < i_0$ and $b_i = 1$ (resp.~$b_i = 2$)
for $i \geq i_0$. In the $D_n$ case, if $b_i = 2$
then one of the balanced cycles in $R_1\cdots R_i$
must be $[1]$.
\end{proposition}

\Proof{Proof.}
By example~\ref{ex:comp}, if the  multiplication of 
$\alpha \in C_n$ by $R \in \mathcal{R}_c$
increases the number of balanced cycles then $l(\alpha R) = l(\alpha)+1$
and $\alpha R$ contains either $1$ or $2$ balanced cycles 
more than $\alpha$. Conversely, if multiplication of $\alpha$
by $R$ decreases either the number of balanced cycles
or the size of a balanced cycle,
 then $l(\alpha R) = l(\alpha)-1$. 
 Since $l(R_1\cdots R_i) + 1 = l(R_1 \cdots R_{i+1})$
 it follows that $b_{i+1}-b_{i} \in \{0,1,2\}$.
As $\gamma_C$ consists of a single balanced cycle, the claim
for $C_n$ is immediate.
For $\gamma_D$, none of the $R_i$ can be of the form~$[j]$
and hence $b_{i+1}-b_{i}$ cannot be $1$.  As the passage
from $R_1\cdots R_i$ to $R_1 \cdots R_{i+1}$ cannot
decrease the size of any balanced cycle and as
$\gamma_D$ contains the balanced cycle~$[1]$, this cycle
must be present in $R_1 \cdots R_{i}$ for each $i \geq i_0$.
\hfill $\Eop$

\begin{corollary}
\label{precedes}
If $\alpha \le \gamma_C$ in $C_n$
then $\alpha$ has at most one balanced cycle.
If $\beta \le \gamma_D$ in $D_n$ then $\beta$
has either no balanced cycles or two balanced
cycles. In the latter case, one of these balanced cycles is $[1]$.
\end{corollary}
\subsection{The $C_n$ lattice.}
Set $\gamma = \gamma_C = [1,2,\dots, n]$.

\begin{definition}
The action of $\gamma$  defines a cyclic order on the set 
$ A=\{1,\ldots,n,-1,\ldots,-n\}$ in which the successor of $i$
is $\gamma(i)$ (thus $1$ is the successor of $-n$). An ordered set
of elements $i_1,i_2,\ldots,i_s$ in $A$ is oriented consistently
(with the cyclic order on $A$)  if there exist integers
$0<r_2<\ldots<r_s \leq 2n-1$ such that $i_j = \gamma^{r_j}(i_1)$
for $j = 2,\ldots,s$. A cycle $\prd{i_1,\ldots,i_s}$
or $[i_1,\ldots,i_s]$ is oriented consistently if the ordered set
$i_1,\ldots,i_s,-i_1,\ldots,-i_s$ in $A$ is oriented consistently.
\end{definition}

\begin{definition}
Two disjoint reflections~$R_1 = \prd{i,j}$ and 
$R_2 = \prd{k,l}$ (resp.~$R_2 = [k]$)
are said to cross if one of the following four ordered sets is oriented
consistently in $A$:
$i , k ,j , l$ or $i , -k , j , -l$ or
$k , i,l,j$ or $k , -i,l,-j$ (resp. $i , k , j ,-k$ or $i , -k ,j , k$
or $k,i,-k,j$ or $-k,i,k,j$).
Two disjoint cycles~$\zeta_1$ and $\zeta_2$ in $C_n$
are said to cross if there exist crossing reflections~$R_1$ and $R_2$
which are 
contained in $\zeta_1$ and $\zeta_2$ respectively.
An element $\sigma \in C_n$
is called \emph{crossing} if some pair of disjoint cycles of $\sigma$
cross. Otherwise $\sigma$ is \emph{non-crossing}.
\end{definition}

 \begin{proposition}\label{le1cn}
 If $\sigma \in C_n$ satisfies $\sigma \le \gamma$ 
 then the cycles of $\sigma$ are oriented consistently  and are 
 noncrossing.
 \end{proposition}

\Proof{Proof.}
We will proceed by induction on $n -l(\sigma)$. If $l(\sigma) = n$ 
then $\sigma = \gamma$ 
and the two conditions of the conclusion are satisfied.  

We assume therefore that the proposition is true for $\tau \in C_n$ 
with 
$n-l(\tau) = 0, 1, \dots , k-1$ and that $\sigma \le \gamma$ satisfies 
$l(\sigma) = n-k$.  By 
definition there is an expression for $\gamma$ as a product of $n$ 
reflections 
$ \gamma = R_1R_2 \dots R_{n-k} R R_{n-k+2} \dots R_n$ with 
$\sigma = R_1R_2\dots R_{n-k}$.  We 
define $\tau = \sigma R$  so that $l(\tau) =  l(\sigma)+1$ and $\tau \le 
\gamma$.  By induction, 
the cycles of $\tau$ are noncrossing and oriented consistently with 
$\gamma$.  

We know that $R$ is either of the form  $ \prd{ i, j  }$ or 
$[i]$ and that $R \leq \tau \leq \gamma$.  Lemma~\ref{contained}
thus implies that $R$ is contained in some paired cycle or some
balanced cycle of $\tau$. The effect of multiplying this cycle
by $R$ is thus described by one of
the first  three equations in Example~\ref{ex:comp}.  Since 
the cycles of $\tau$ 
are noncrossing and oriented consistently with $\gamma$, we 
see that the same is true for $\sigma$.
\hfill $\Eop$

\begin{proposition}\label{le2cn}
Let   $\sigma \in C_n$. If the cycles of $\sigma $ are 
oriented consistently and are 
noncrossing then $\sigma \leq \gamma$.
\end{proposition}

\Proof{Proof.} 
Assume that $\sigma \in C_n$ satisfies 
the two hypotheses of the proposition.
Write $\sigma = c_1\bar{c}_1 \ldots c_a \bar{c}_a \gamma_1 \ldots \gamma_b$
and set $t(\sigma) = a+b$. We proceed by induction on $t(\sigma)$.
If $t(\sigma)=1$ then either $\sigma$ consists
of a single balanced cycle or a single paired cycle.
In the former case, consistent orientation implies that 
$\sigma = \gamma$. In the latter case, consistent orientation
implies that 
$\sigma = \prd{ i,i+1,\ldots,n,-1,\ldots,-i+1  }$ for some $i$.
As $l(\sigma)=n-1$ and 
$\sigma [i-1] = \gamma$, we see that $\sigma \leq \gamma$. 
\vskip .2cm
Assume now that $t(\sigma) \geq 2$ and that the proposition is
true for each element~$\theta\in C_n$ with $t(\theta) < t(\sigma)$.
If $\sigma$ contains a balanced cycle, the non-crossing
hypothesis implies that there can be only one which we denote
$\tau = [i_1,\ldots,i_r]$. Otherwise let 
$\tau = \prd{i_1,\ldots,i_r}$ be some paired cycle of $\sigma$.
As $\sigma \neq \tau$, there exists an $i_k$ whose successor
does not lie in $\{\pm i_1,\ldots, \pm i_r\}$. By choosing one
 of the other $2r-1$ cycle expressions for $\tau$ if necessary,
we may assume that the successor~$j_1$ of $i_r$
does not lie in $\{\pm i_1,\ldots, \pm i_r\}$. Let
$\rho = \prd{j_1,\ldots,j_s}$ be the paired cycle of $\sigma$
which contains $j_1$ and let $R=\prd{i_r,j_s}$.
Then $\sigma = \tau \rho \sigma_1 \ldots \sigma_k$ for some disjoint
paired cycles $\sigma_1,\ldots,\sigma_k$ (some $k \geq 0$)
and
\[ \sigma R = \left\{ \begin{array}{ll}
 [i_1,\ldots,i_r,j_1,\ldots,j_s]\sigma_1 \ldots \sigma_k & \mbox{ or} \\
\prd{i_1,\ldots,i_r,j_1,\ldots,j_s}\, \sigma_1 \ldots \sigma_k. & 
\end{array} \right. \]
Note that $t(\sigma R) = t(\sigma) -1$. As the cycles $\tau$ and $\rho$
do not cross and each is oriented consistently, our choice
of $j_1$ ensures that the ordered set $i_1,\ldots,i_r,j_1,\ldots,j_s, 
-i_1,\ldots,-i_r,-j_1,\ldots,-j_s$ is also
oriented consistently.

Assume now that one of the cycles $\sigma_e$ crosses the cycle
$\tau \rho R$  of $\sigma R$. Then there exist crossing reflections
$R_1$ and $R_2$ contained in $\tau \rho R$ and $\sigma_e$
respectively. Since $\sigma_e$ is paired, $R_2$ is necessarily
paired; $R_2 = \prd{c,d}$ say.
  Since $\sigma$ is non-crossing, $R_1$ cannot
be contained in $\tau$ or in $\rho$. There are three cases to consider
\begin{enumerate}
\item $R_1 = \prd{i_a,j_b}$
for some $1 \leq a \leq r$ and $1 \leq b \leq s$.
\item $R_1 = \prd{j_b,-j_b} $ for some $1 \leq b \leq s$ ($\tau$ is necesarily balanced ).
\item $R_1 = \prd{i_a,-j_b} $ for some $1 \leq b \leq s$ ($\tau$ is necesarily balanced ).
\end{enumerate}
By a suitable choice of the representative 
$R = \prd{c,d} = \prd{d,c}=\prd{-c,-d} =\prd{-d,-c}$, the first case splits into
two essential subcases: {\it (a)} the ordered set $i_a,c,j_b,d$ is oriented consistently
and {\it (b)} the ordered set $c,i_a,d,j_b$ is oriented consistently.
 We know that $c$ is not in
$\{ \pm i_1,\ldots,\pm i_r , \pm j_1, \ldots, \pm j_s\}$.
In particular $c \neq i_r,\ j_1$. In case $(1a)$, if $c$   precedes
$i_r$, then $S= \prd{i_1,i_r}$  is contained in $\tau$ and crosses
$R_2$, contradicting the fact that $\sigma$ is noncrossing.
Likewise, if $c$ follows $i_r$ then $c$ follows $j_1$ and $S = \prd{j_1,j_b}$
is contained in $\rho$ and crosses $R_2$, again contradicting
the fact that $\sigma$ is non-crossing. Thus case $(1a)$ is impossible.
A similar argument shows that case $(1b)$ is also impossible.

As in case $1$, case $2$ splits into two subcases: {\it (a)} the ordered set
$j_b,c,-j_b,d$ is oriented consistently and {\it (b)} the ordered set
$c,j_b,d,-j_b$ is oriented consistently. In case $(2a)$, if $c$ precedes $-i_r$ then
 the ordered set
$ i_r , j_b , c, -i_r,d $ is oriented consistently
and hence $\prd{c,d}$ crosses $[-i_r] \subset \tau$.
But this contradicts the fact that $\sigma$ is non-crossing.
If $c$ follows $-i_r$, then $c$ necessarily succeeds $-j_1$
and we find that the ordered set $-j_1,c,-j_b , d$ is consistently
oriented. Thus $\prd{c,d}$
crosses $\prd{-j_1,-j_b} \subset \rho$, again contradicting
the fact that $\sigma$ is non-crossing. Thus case $(2a)$ is impossible.
A similar argument shows that case $(2b)$ is also impossible.

Finally, case 3 also splits into two subcases:
{\it (a)} the ordered set
$i_a,c,-j_b,d$ is oriented consistently and {\it (b)} the ordered set
$c,i_a,d,-j_b$ is oriented consistently.
We show that $(3b)$ is impossible (the proof that case $(3a)$ is impossible
is similar). We are given that the ordered set
$c,i_a,d,-j_b$ is oriented consistently. If $d$ precedes $-i_a$ then
$\prd{c,d}$ crosses $[i_a]$ in $\sigma$, a contradiction.
Therefore $d$ follows $-i_a$. If $d$ now precedes $-i_r$, then the ordered set
$c,-i_a,d,-i_r$ is oriented consistently. Hence $\prd{-i_a,-i_r}$
crosses $\prd{c,d}$ in $\sigma$, a contradiction. Therefore $d$
follows $-i_r$ and hence $-j_1$. But now $\prd{-j_1,-j_b}$  crosses
$\prd{c,d}$ in $\sigma$, a contradiction. Thus case $(3b)$ is impossible.

We conclude that the cycles $\tau \rho R$ and $\sigma_e$ do not cross.
Since no two distinct elements of $\sigma_1,\ldots,\sigma_k$ cross
(because $\sigma$ is assumed non-crossing),
it follows that $\sigma R$ is non-crossing. As
$t(\sigma R) = t(\sigma) -1$ and the cycles of $\sigma R$
are oriented consistently, it follows by induction
that $\sigma R \leq \gamma$. Thus there exist reflections $R_1,
\ldots, R_k$ with $k = n -l(\sigma R)$ and
\begin{equation}
\label{star}
\sigma R R_1 \ldots R_k = \gamma
\end{equation}
As $l(\sigma R) = l(\sigma) +1$ by lemmas \ref{pairedlength} and \ref{balancedlength} 
and proposition \ref{dimfix}, 
we see that $k+1 = n - l(\sigma)$.  Hence equation~(\ref{star})
also implies that $\sigma \leq \gamma$. \hfill $\Eop$

\begin{lemma}
If $\sigma \leq \gamma$ and $\tau \leq \gamma$ then
$\sigma \leq \tau$ if and only if
$\sigma \subset \tau$.
\end{lemma}

\Proof{Proof.}  Follows from Lemma~\ref{orderchar} and lemma~\ref{contained}.
\hfill 
$\Eop$
\vskip .2cm
Combining the previous three results yields the following Theorem.
\begin{thm}
Let $NCP$ denote Reiner's non-crossing partition lattice
for the $C_n$ group from \cite{reiner}. The mapping
\[ \colon \{\alpha \in C_n : \alpha \leq \gamma \} \longrightarrow NCP \]
which takes $\alpha$  to the noncrossing partition defined by its 
cycle structure is a bijective poset map. In particular,
$\{\alpha \in C_n : \alpha \leq \gamma \}$ is a lattice.
\end{thm}
\subsection{The $D_n$ lattice.}
Set $\gamma = \gamma_D = [1][2,3, \dots, n]$ and suppose $\alpha \le \gamma$.  
Recall from Corollary~\ref{precedes} that for such an $\alpha$ either 
$[1][k] \le \alpha$ for some $k \in \{2,3,\dots , n\}$ or $l$ and $-l$ are in 
different $\alpha$ orbits for all $l \in \{1,2,\dots , n\}$.  In the former case 
we will call $\alpha$ {\em balanced} and in the latter case we will call $\alpha$ 
{\em paired}.
\vskip .2cm
We note that lattices are associated to the groups $C_n$ and $D_n$ in \cite{reiner}.  
We have shown the Reiner $C_n$ lattices are isomorphic to ours.  However the Reiner 
$D_n$ lattices are not the same as the ones we consider.  In particular, the Reiner 
$D_n$ lattices are subposets of the Reiner $C_n$ lattices. 
\vskip .2cm
To show that the interval $[I, \gamma]$ in $D_n$ is a lattice we will compute 
$\alpha \wedge \beta$ for $\alpha , \beta \le \gamma$.  Since the poset is finite 
the existence of  least upper bounds follows.  We will 
consider different cases depending on the types of $\alpha$ and $\beta$.  In all cases we 
will construct a candidate $\sigma$ for $\alpha \wedge \beta$ and show that 
$\sigma \in D_n$, $\sigma \le \alpha, \beta$ and 
$S_{\alpha} \cap S_{\beta} \subset S_{\sigma}$.  
Since the reverse inclusion is immediate it follows from Lemma~\ref{meetchar} that 
$\sigma = \alpha \wedge \beta$.
\vskip .2cm
\begin{note}
In this section we will frequently pass between the posets determined by $C_n$, $D_n$ and 
several other finite reflection subgroups of $C_n$.  As the partial order on each of 
these groups is the restriction of the partial order on $O(n)$, we can use the same 
symbol $\le $ to denote the partial order in each case.  The reflection subgroup 
in question should be clear from the context.
\end{note}
\vskip .2cm
Suppose first that both $\alpha$ and $\beta$ are balanced. 
Since $D_n \subset C_n$ and $C_{n-1}$ can be identified with the subgroup of $C_n$ which fixes 
$1$, each balanced element of $D_n$ can be used to define a balanced element of $C_{n-1}$, 
that is, an element containing a balanced cycle.  
Thus we define the balanced $C_{n-1} $ elements $\alpha'$ and $\beta'$ by
\[\alpha = [1]\alpha' \quad \mbox{ and } \quad \beta = [1]\beta'\] 
and the $C_{n-1}$ element $\sigma' = \alpha' \wedge \beta'$, where the meet is 
taken in $C_{n-1}$.  Now $\sigma'$ may or may not be balanced.  If $\sigma' $ is balanced 
define the $C_n$ element $\sigma$ by $\sigma = [1]\sigma'$.  If $\sigma'$ is not balanced 
set $\sigma = \sigma'$.  
\begin{proposition} If $\alpha$ and $\beta $ are balanced and $\sigma$ is defined as above 
then $\sigma \in D_n$, $\sigma \le \alpha, \beta $ and 
$S_{\alpha} \cap S_{\beta} \subset S_{\sigma}$.
\end{proposition}
\Proof{Proof.}
We show that $\sigma \in D_n$ and $\sigma \le \alpha$.   The proof that 
$\sigma \le \beta$ is completely analogous.  
First consider the case 
where $\sigma'$ is balanced.  Thus $[k] \le \sigma' \le \alpha'$ in $C_{n-1}$ 
for some $k $ satisfying $2 \le k \le n$.  
So we can find reflections $R_1, \dots , R_s$ in $C_{n-1}$ with
\[\alpha' = R_1R_2\dots R_s, \quad \quad \sigma' =  R_1R_2\dots R_t, \quad \quad 
 R_1 = [k],\]
where $l(\alpha') = s \ge t = l(\sigma')$.  Since $\alpha' \in C_{n-1}$, 
Lemma~\ref{balancedlength} gives $R_{2}, \dots , R_s$ all of the form 
$\prd{i,j}$ or $\prd{i, -j}$ for $2 \le i < j \le n$.  In particular, these reflections 
lie in $D_n$.  Now $\alpha$ is of length $s+1$ in $C_n$ and 
\begin{eqnarray*}
\alpha  &=& [1]R_1R_2\dots R_tR_{t+1}\dots R_s\\
&=& [1][k]R_2\dots R_tR_{t+1}\dots R_s\\
&=& \prd{1,k}\prd{1,- k}R_2\dots R_tR_{t+1}\dots R_s.
\end{eqnarray*}
This last expression only uses $D_n$ reflections so that 
\[\sigma = \prd{1,k}\prd{1, - k}R_2\dots R_t \le \alpha \quad \mbox{ in } D_n.\] 
\vskip .2cm
Next we consider 
the case where $\sigma'$ is paired.  Here $\sigma' \le \alpha'$ and $\alpha'$ is balanced so 
we can find reflections $R_1, \dots , R_s$ in $C_{n-1}$ with
\[\alpha' = R_1R_2\dots R_s, \quad \quad \sigma' =  R_1R_2\dots R_t, \]
where $l(\alpha') = s > t = l(\sigma')$ and exactly one of 
$R_{t+1}, \dots, R_s$ is of form $[k]$.  Since $R[k] = [k]([k]R[k])$, we can 
assume $R_{t+1} = [k]$.  Note also that $R_1, \dots , R_t$ are 
each of the form $\prd{i,j}$ or $\prd{i, -j}$ for $2 \le i < j \le n$ 
and hence commute with $[1]$ in $C_n$.  Thus we can write the 
following identities in $C_n$. 
\begin{eqnarray*}
\alpha  &=& [1]R_1R_2\dots R_t[k]R_{t+2}\dots R_s\\
&=& R_1\dots R_t[1][k]R_{t+2}\dots R_s\\
&=& R_1\dots R_t\prd{1,k}\prd{1,- k}R_{t+2}\dots R_s.
\end{eqnarray*}
This last expression only uses $D_n$ reflections so that $\sigma \le \alpha $ in $D_n$.  
\vskip .2cm
Finally we show that 
$S_{\alpha}\cap S_{\beta}\subset S_{\sigma}$.  First suppose $\sigma'$ is balanced 
and $R \in S_{\alpha}\cap S_{\beta}$.  Thus $R$ is a 
reflection satisfying $R\le \alpha , \beta$.  
If $R$ is of the form $\prd{1,k}$, then $[1][k] \le \alpha, \beta$ since $k$ must belong 
to a balanced cycle of both $\alpha$ and $\beta$.  
Thus $[k] \le \alpha', \beta'$ so that $[k] \le \sigma'$ 
and $[1][k] \le \sigma$, which gives $\prd{1,k} \le \sigma$ as required.  If $R$ 
is not of form $\prd{1,k}$ then $R \le \alpha , \beta$ implies $R \le \alpha' , \beta'$ so 
that $R \le \sigma'$ and $R \le \sigma$.  
\vskip .2cm
In the case where $\sigma'$ is paired, 
$R \le \alpha, \beta$ implies $R$ must be of form $\prd{i,j}$ or $\prd{i, -j}$ 
for $2 \le i < j \le n$ so that $R \le \alpha', \beta'$ giving $R \le \sigma' = \sigma$.
\ \hfill \Eop
\vskip .2cm
Since we have completed the case where both $\alpha$ and $\beta$ are balanced we will 
assume from now on that $\alpha$ is paired.  We note some consequences of this fact 
which will apply in the remaining cases.  The fact that $\alpha$ is paired means that 
$\alpha \le \prd{1,k}\gamma$ or $\alpha \le \prd{1, -k}\gamma$ for some $k \in 
\{2, 3, \dots, n\}$.  
Since conjugation by the $C_{n-1}$ element $[2, \dots, n]$ is a poset isomorphism of 
the interval $[I, \gamma]$ in $D_n$, 
we may assume for convenience of notation that $k = - 2$ so that 
\[\alpha \le \prd{1, - 2}[1][2, \dots , n] = \prd{1,2,\dots ,n}.\]
If we let $\delta = \prd{1,2,\dots ,n}$ then a reflection $R$ in $D_n$ satisfies 
$R \le \delta $ if and only if $R \subset \delta$.   Thus we can 
identify the interval $[I,\delta]$ in $D_n$ with the set of non-crossing partitions of 
$\{1,2, \dots , n\}$.  Recall that a non-crossing partition of the ordered set 
$\{a_1,a_2, \dots , a_n\}$ is a partition with the property that whenever
\[1 \le i < j < k < l \le n\]
with $a_i,a_k$ belonging  to the same block $B_1$ and 
$a_j,a_l$ belonging  to the same block $B_2$ we have $B_1 = B_2$. 
If $\alpha \wedge \beta$ exists, it 
will satisfy 
\[\alpha \wedge \beta \le \alpha \le \prd{1,2,\dots , n}\] 
and so will correspond 
to a noncrossing partition of $\{1,2,\dots, n\}$.  Accordingly, we define a 
reflexive, symmetric relation on $\{1,2,\dots, n\}$ by 
\[i \sim j \quad \Leftrightarrow \quad i = j \quad \mbox{ or } \quad 
\prd{i,j} \le \alpha , \beta. \]
We need to show that $\sim$ is transitive and hence is an equivalence relation.
We then show that the resulting partition of $\{1,2,\dots, n\}$ is noncrossing and 
 determines an element $\sigma$ of $D_n$ which satisfies $\sigma \le \alpha, \beta$ and 
 $S_{\alpha}\cap S_{\beta} \subset S_{\sigma}$.
\vskip .2cm
Suppose that $\alpha$ is paired and $\beta$ is balanced.  
Recall that $\beta$ has two balanced cycles, 
one of which is $[1]$.  For convenience of terminology we will call the other balanced 
cycle the second balanced cycle of $\beta$.   As above we will 
have occasion to use the balanced element~$\beta' \le [2, \dots , n]$ in $C_{n-1}$ defined by 
$\beta = [1]\beta'$. 
\begin{proposition} If $\alpha$ is paired and $\beta $ is balanced then 
the relation $\sim$ above determines an element $\sigma$ of $D_n$ satisfying 
 $\sigma \le \alpha, \beta$ and $S_{\alpha}\cap S_{\beta} \subset S_{\sigma}$.
\end{proposition}
\Proof{Proof.}
First we establish the transitivity of the 
$\sim$ relation.  Suppose $i,j,k$ are distinct elements of $\{1,2,\dots, n\}$ with 
$i\sim j$ and $j \sim k$.  Since $\prd{i,j}, \prd{j,k} \le \alpha$ we get 
$\prd{i,k} \le \alpha$ since $\alpha$ corresponds to a partition of $\{1,2,\dots, n\}$.  
If $1 \not \in \{i,j,k\}$ then $\prd{i,j}, \prd{j,k} \subset \beta$ 
(s-containment cannot arise) and it follows that $\prd{i,k} \le \beta$.  
If $i = 1$, then $\prd{i,j} \le \beta$ means that 
$\prd{i,j} \sqsubset \beta$ so that 
$j$ belongs 
to the second balanced cycle of $\beta$.  Since $j\sim k \ne 1$, $k$ also belongs to this second 
balanced cycle and $\prd{i,k}\le [1][j,k] \le \beta$.  If $j = 1$, then both $i$ and $k$ belong 
to the second balanced cycle of $\beta$.  Hence 
$\prd{i,k} \le \beta$.  The case $k=1$ is analogous to the case $i = 1$.  
\vskip .2cm
To show that the partition of $\{1, \dots , n\}$ defined by $\sim$ is non-crossing suppose 
$1 \le i < j < k < l \le n$ with 
\[\prd{i,k}, \prd{j,l} \le \alpha, \beta.\]
Since $\alpha$ corresponds to a noncrossing partition we have $\prd{i,j,k,l} \le \alpha$.  
If $i = 1$, then $k$ belongs to the second balanced cycle and $[k] \le \beta'$ in 
$C_{n-1}$.  Since 
$1 < j < k < l $, $\prd{j,l} \le \beta'$ and $\beta' \le [2,\dots, n]$ in $C_{n-1}$, the 
crossing pair consisting of $(j,l)$ and $(k, -k)$ must lie in the same $\beta'$ cycle.  
Thus $[j,k,l] \le \beta'$ and $\prd{1,j,k,l} \le [1][j,k,l] \le \beta$.  If $i \ne 1$, then 
$\prd{i,k}, \prd{j,l} \le \beta'$ and since $\beta' \le [2,\dots, n]$ in $C_{n-1}$, 
$\prd{i,j,k,l} \le \beta'$ by 
proposition~\ref{le1cn}, giving $\prd{i,j,k,l} \le \beta$.
\vskip .2cm
Thus the relation $\sim$ defines a noncrossing partition of $\{1,2,\dots , n\}$ and hence 
determines an element $\sigma$ of $D_n$.  
By the definition of $\sim$ the element $\sigma$ satisfies  $\sigma\le \alpha, \beta$ and 
$S_{\alpha}\cap S_{\beta} \subset S_{\sigma}$. \hfill \Eop
\vskip .2cm
Finally we consider the case where both $\alpha$ and $\beta$ are paired. 
\begin{proposition} If $\alpha$ and $\beta $ are paired then 
the relation $\sim$ above determines an element $\sigma$ of $D_n$ satisfying   
$\sigma \le \alpha, \beta$ and  $S_{\alpha}\cap S_{\beta} \subset S_{\sigma}$.
\end{proposition}
\Proof{Proof.}
To establish the transitivity of $\sim$ in this case 
let $i,j,k$ be distinct elements of $\{1,2,\dots, n\}$ with 
$i\sim j$ and $j \sim k$.  As in the previous proposition, 
$\prd{i,k} \le \alpha$ follows immediately.  
Since $\beta$ is paired, $i\sim j$ and $j \sim k$ mean that $i,j,k$ belong to the same 
cycle of $\beta$ so that $\prd{i,k} \le \beta$ also.
\vskip .2cm
To show that the partition of $\{1, \dots, n\}$ defined by $\sim$ is noncrossing suppose 
$1 \le i < j < k < l \le n$ with 
\[\prd{i,k}, \prd{j,l} \le \alpha, \beta.\]
Since $\alpha$ corresponds to a noncrossing partition we have $\prd{i,j,k,l} \le \alpha$.  
The element $\beta$ is paired so 
we can assume $\beta \le \tau = \prd{1,m}\gamma$ or $\beta \le \tau = \prd{1,-m}\gamma$, for 
some $m \in \{2, 3, \dots , n\}$.  Looking at the case $\tau = \prd{1,m}\gamma$ first we 
get 
\[ \tau =    \prd{1, -m, -m-1, \dots, -n, 2, 3, \dots ,m-1}. \]
Since $\beta \le \tau$ the element $\beta$ corresponds to a noncrossing partition 
of the ordered set $\{1, -m, -m-1, \dots, -n, 2, 3, \dots ,m-1\}$.  
Since $1 \le i < j < k < l \le n$, we deduce that either 
\[1 \le i < j < k < l \le m-1 \quad \mbox{ or } \quad m \le i < j < k < l \le n.\]
Since $\beta $ corresponds to a noncrossing partition 
of the ordered set 
\[\{1, -m, -m-1, \dots, -n, 2, 3, \dots ,m-1\}\]
 and 
$\prd{i,k}, \prd{j,l} \le  \beta$ it follows in either case that 
$\prd{i,j,k,l} \le \beta$. 
 The case $\tau = \prd{1,-m}\gamma$ is similar.  Here
\[ \tau =    \prd{1, m, m+1, \dots, n, -2, -3, \dots ,-m+1}, \]
and again we can deduce $\prd{i,j,k,l} \le \beta$.
\vskip .2cm
Thus $\sim$ defines a noncrossing partition of $\{1,2,\dots , n\}$ and hence 
an element $\sigma$ in $D_n$ satisfying  $\sigma \le \alpha, \beta$ 
and $S_{\alpha}\cap S_{\beta} \subset S_{\sigma}$ 
as in previous proposition. \hfill \Eop
\vskip .2cm
Combining the results of this subsection we obtain the following theorem.
\begin{thm}
The interval $[I, \gamma]$ in $D_n$ is a lattice.
\end{thm}
\section{ Poset groups and K$(\pi,1)$'s.}

\begin{defn}
If $W$ is a finite Coxeter group  and $\gamma \in W$ we define the {\em poset group} 
$\Gamma = \Gamma(W, \gamma)$ to be the group with the following presentation.  
The generating set for $\Gamma$ consists of a 
copy of the set of non-identity elements in $[I,\gamma]$.  We will denote by 
$\{w\}$ the generator of $\Gamma$ corresponding the element $w \in (I, \gamma]$.  
The relations in $\Gamma$ are  all identities of the form 
$\{w_1\}\{w_2\} = \{w_3\}$, where
$w_1, w_2$ and $ w_3 $ lie in $(I, \gamma]$ with $w_1 \le w_3$ and 
$w_2 = w_1^{-1}w_3$.
\end{defn}
Since none of the relations involve inverses of the generators, there is a semigroup, 
which we will denote by $\Gamma_+ = \Gamma_+(W, \gamma)$, with the same presentation.  
As in section $5$ of \cite{brady1}, we define a \emph{positive} word in $\Gamma$ 
to be a word in the generators that does not involve the inverses of the
generators.  We say two positive words $A$ and $B$ are 
\emph{positively equal}, and we write $A \doteq B$, if $A$ can be transformed to $B$ through a
sequence of positive words, where each  word in the sequence is obtained
from the previous one by replacing one side of a defining relator by the other side.  
Since the interval $(I, \gamma]$ inherits the reflection length from $W$ we use this to 
associate a \emph{length} to each generator of $\Gamma(W, \gamma)$ and hence a length $l(A)$ to 
each positive word $A$.  It is immediate that positively equal words have the same length. 
\vskip .2cm
From now on we only consider those pairs $(W,\gamma)$ with the property that  
\emph{ the interval $[I,\gamma]$ in $W$ forms a lattice}.  
It is clear that the results stated for  the braid group 
in sections $5$ and $6$ of \cite{brady1} apply to poset groups under this extra assumption.  
We will review them briefly below. 
\vskip .2cm
In \cite{brady1} it is shown that this lattice condition is satisfied when $W$ is a 
Coxeter group of type $A_n$ and $\gamma$ is a Coxeter element.  In section 4 above we 
have shown that the lattice condition is satisfied when $W$ is a 
Coxeter group of type $C_n$ or $D_n$  and $\gamma$ is a Coxeter element.   When the 
Coxeter group is generated by two reflections the lattice condition is automatic for 
any $\gamma$.  When the Coxeter group is generated by three reflections the lattice 
condition reduces to checking the only case where 
\[\alpha \wedge \beta  \not\in \{\alpha, \beta, \gamma\}.\]
This occurs when $\alpha$ and $\beta$ are distinct reflections and have at least one 
common upper bound of length $2$.  Any such length $2$ element $\delta$ must have 
$F(\delta)$ coinciding with  the unique line 
of intersection of the two reflection planes.  Hence $\delta$ is unique.  
This is precisely the ingredient which makes the metric constructed in \cite{brady} 
have non-positive curvature.

The following result is taken from \cite{brady1}.  Its proof is the same.
\begin{lemma}  Assume that the interval  $[I, \gamma]$ 
forms a lattice and suppose $a,b, c \le \gamma$.  We define nine elements 
$d$, $e$, $f$, $g$, $h$, $k$, $l$, $m$ and $n$ in $[I, \gamma]$ 
 by the equations 
\[ a\vee b = ad = be, \quad b\vee c = bf = cg,
\quad c\vee a = ch = ak\]
and
\[a\vee b\vee c = (a\vee b)l = (b\vee c)m = (c\vee a)n. \]
Then we can deduce
\[e\vee f = el = fm, \quad d\vee k = dl = kn, \quad h\vee g = hn = gm.  \]
\end{lemma}
The statements and proofs of the results of section~5 and section~6 of \cite{brady1} 
generalize in a straightforward manner to the current setting.  In particular, we have 
the following definitions and results.
\begin{lemma}
The semigroup associated to $\Gamma$ has right and
left cancellation properties.
\end{lemma}
\begin{lemma}
Suppose $a_1,a_2, \dots ,a_ k \le \gamma$ in $W$,
$P$ is positive and
\[ P \doteq X_1\{a_1\} \doteq \dots \doteq X_ k\{a_k\} \]
with $X_ i$ all  positive.  Then there is a positive word $Z$ satisfying
\[P \doteq Z\{a_1\vee \dots \vee a_ k\}.\]
\end{lemma}
\begin{thm}
In $\Gamma$, if two positive words are equal they
are positively equal.  In other words, the semigroup $\Gamma_+$ embeds in $\Gamma$.
\end{thm}
As in \cite{brady1} we define an abstract simplicial complex $X(W, \gamma)$ for each  
$\Gamma(W, \gamma)$.
\begin{defn}
We let $X = X(W, \gamma)$ be the abstract simplicial complex 
with vertex set $\Gamma$, which has a $k$-simplex on the subset 
$\{g_0, g_1, \dots , g_k\}$ if and only if $g_i = g_0\{w_i\}$ for 
$i = 1, 2, \dots , k$ where 
\[I < w_1 < \dots < w_k \le \gamma \quad \mbox{ in $W$}.\]
\end{defn}
There is an obvious simplicial action of $\Gamma$ on $X$ given by 
\[ g\cdot\{g_0, g_1, \dots , g_k\} = \{gg_0, gg_1, \dots, gg_k\}. \]
The main result of section $6$ of \cite{brady1} also holds for these poset groups.
\begin{thm}
$X(W, \gamma)$ is contractible. 
\end{thm}
If we define $K = K(W, \gamma)$ to be the quotient space $K = \Gamma \backslash X$, then $K$ is  
a $K(\Gamma, 1)$. 
\vskip .2cm
We finish this section with an example of a poset group $\Gamma(W, \gamma)$, 
with $[I, \gamma]$ a lattice but $\gamma$ not a Coxeter element in $W$.
\begin{example}  Let $W = C_2$ and $\gamma = [1][2]$.  The group $\Gamma(C_2, \gamma)$ 
has presentation
\[\langle a, b, c, d, x \mid x = ab  = ba = cd = dc \rangle\]
where $a = \{[1]\}$, $b = \{[2]\}$, $c = \{\prd{1,2}\}$, $d = \{\prd{1,-2}\}$ and  
$x = \{[1][2]\}$.  From the presentation we see that $\Gamma$ is an amalgamated 
free product of a copy $\iIZ \times \iIZ$ generated by $a$ and $b$ with 
a copy $\iIZ \times \iIZ$ generated by $c$ and $d$ over the infinite cyclic subgroup 
generated by $x$.  The above construction gives a two-dimensional contractible 
universal cover for the presentation $2$-complex which can be shown to be simplicially 
isomorphic to $X(C_2, [1,2])$.
\end{example}
\section{Group Presentations.}  
In this section we prove that the poset groups $\Gamma(W,\gamma)$ of section 5 are 
isomorphic to the Artin groups $A(W)$ for $W $ of type $C_n$ or $D_n$ and $\gamma$ the 
appropriate Coxeter element.  
The proof is based on the following surprising 
property that these Artin groups share with the braid group.  If $X = x_1x_2\dots x_n$ 
is the product of the standard Artin generators then there is a finite set of elements 
in $A(W)$ which is invariant under conjugation by $X$.  Moreover under the canonical 
surjection from $A(W)$ to $W$ this set is taken bijectively to the set of reflections 
in $W$. 
The following lemma is a straightforward generalisation of Lemma~4.5 of \cite{brady1}.
\begin{lem}\label{newpres}
The poset group $\Gamma(W, \gamma)$ 
is isomorphic to the abstract group generated by the set of all $\{R\}$, for $R$ a reflection 
in $[I, \gamma]$, 
subject to the relations 
\[
\{R_1\}\{R_2\}\dots \{R_n\} = \{S_1\}\{S_2\}\dots \{S_n\}, 
\] 
for $R_i, S_j$ reflections satisfying 
\[\gamma = R_1R_2\dots R_{n} \quad  \mbox{ and } \quad 
\gamma = S_1S_2\dots S_{n},\]
where $n = l(\gamma)$.  
\end{lem}
We will refer to $\{w\}\in \Gamma(W, \gamma)$ as the lift of $w \in W$ whenever $w \le \gamma$.  
In particular, we  will refer to $\{w\}$ as a reflection lift whenever $w $ is 
a reflection.
\vskip .2cm
Since the Artin groups of type $C_n$ and $D_n$ both contain copies of the $n$-strand braid 
group $B_n$ we collect here some facts about the braid group which will be useful.  We recall 
that $B_n$ is the group with generating set $x_2, x_3, \dots x_n$ and defining relations 
\[x_ix_{i+1}x_i = x_{i+1}x_ix_{i+1}\quad \mbox{ for } \quad 2 \le i \le n-1,\]
\[x_ix_j = x_jx_i \quad \mbox{ for } \quad |j-i| \ge 2.\]
We define $x_{i,j}$ and $Y_{i,j}$, for $1\le i < j \le n$  by 
\[Y_{i,j} = x_{i+1}\dots x_j, \quad \mbox{ and } \quad Y_{i,j} = Y_{i+1,j}x_{i,j}.\]
Then Lemma~4.2 of \cite{brady1} gives, for $1 \le i < j<k \le n$,
\[x_{i,j}x_{j,k} = x_{j,k}x_{i,k} = x_{i,k}x_{i,j}.\]
Since $x_k = x_{k-1,k}$ it follows that $x_{i,j}Y_{i,j-1} = Y_{i,j}$ and that
\[x_kY_{i,j} = Y_{i,j}x_{k-1}\quad \mbox{for} \quad i+2 \le k \le j.\] 
When $k = i+1$ we have $x_{i+1}Y_{i,j} = x_{i+1}Y_{i+1,j}x_{i,j} = Y_{i,j}x_{i,j}$.
\subsection{The $C_n$ case}
The Artin group $A(C_n)$ has a presentation with generating set 
$x_1, x_2, \dots x_n$, 
subject to the relations 
\[x_1x_2x_1x_2 = x_2x_1x_2x_1\]
\[ x_ix_{i+1}x_i = x_{i+1}x_ix_{i+1}\]
whenever $1< i < n$ and 
\[x_ix_j = x_jx_i\]
whenever $|j-i| \ge 2$.  
\begin{defn}  We define a function $\phi$ from the generators of $A(C_n)$ to 
$\Gamma(C_n, \gamma)$ by 
\[x_1 \mapsto \{[1]\}, x_2 \mapsto \{\prd{1,2}\}, x_3 \mapsto \{\prd{2,3}\}, \dots , x_n \mapsto \{\prd{n-1, n}\}\]
\end{defn}
\begin{lem}  The function $\phi$ determines a well-defined and surjective homomorphism.
\end{lem}
\emph{Proof:}  \hspace{.2cm} The relations involving $\phi(x_1)$ hold in $\Gamma(C_n, \gamma)$ by virtue 
of the following identities in $\Gamma(C_n, \gamma)$. 
\begin{eqnarray*}
\{[1]\}\{\prd{1,2}\}\{[1]\}\{\prd{1,2}\} &=& \{[1,2]\}\{[1,2]\} \\
&=& \{\prd{1,2}\}\{[2]\}\{\prd{1, -2}\}\{[1]\}\\
&=& \{\prd{1,2}\}\{[1,2]\}\{[1]\}\\
&=& \{\prd{1,2}\}\{[1]\}\{\prd{1,2}\}\{[1]\}
\end{eqnarray*}
\[\{[1]\}\{\prd{i, i+1}\} = \{\prd{i, i+1}\}\{[1]\}, \quad \mbox{for} \quad i  \ge 2. \]
The image of the subgroup generated by $\{x_2, \dots , x_n\}$ lies in the copy of the braid 
group corresponding to $\Sigma_n < C_n$ so that the relations not involving $\phi(x_1)$ 
hold by Lemma~4.2 and Lemma~4.4 of \cite{brady1}.  Thus $\phi$ is well-defined.

To establish surjectivity, first note that 
\[\{\prd{i, i+1, \dots , j}\}  = \phi(Y_{i,j}) \quad \mbox{ and } \quad 
\{\prd{i,j}\}  = \phi(x_{i,j}) \]
for $1 \le i < j \le n$ all lie in $\textrm{im}(\phi)$.  
Next  $\{[j]\} \in \textrm{im}(\phi)$ since 
\[\phi(x_1x_{1,j} ) = \{[1]\}\{\prd{1,j}\} = \{[1,j]\} = \{\prd{1,j}\}\{[j]\}.\]
Finally, $\{\prd{i,-j}\}\in \textrm{im}(\phi)$ for $1 \le i < j \le n$ since
\[\{\prd{i,j}\}\{[j]\} = \{[i,j]\} = \{[j]\}\{\prd{i,-j}\}.\]
\ \hfill \Eop
\vskip .2cm
To construct an inverse to $\phi$ we will use the presentation for $\Gamma(C_n, \gamma)$ 
given by lemma~\ref{newpres}.
\begin{definition}
We define a function $\theta$ from the generators of $\Gamma(C_n, \gamma)$ to $A(C_n)$ by
\begin{eqnarray*}
\{[1]\} &\mapsto& x_1\, , \,\{\prd{i,j}\} \mapsto x_{i,j}\quad \mbox{ for }\quad 
1 \le i < j \le  n,\\
\{[j]\} &\mapsto& y_j \,\,\mbox{ for }\,\, 2 \le j \le n, \quad 
\{\prd{i,-j}\} \mapsto z_{i,j}\,\, \mbox{ for }\,\, 1 \le i < j \le n,
\end{eqnarray*}
where $y_j$ is the unique element of $A(C_n)$ satisfying
\[x_1x_2\dots x_j  = x_2\dots x_jy_j\]
and $z_{i,j}$ is the unique element of $A(C_n)$ satisfying
\[  z_{i,j}y_i = y_ix_{i,j}. \]
\end{definition}
The homomorphism determined by  $\theta$ will be  surjective since each $x_i$ is the 
image of some reflection lift. 
We note that $Y_{i,j}y_j = y_iY_{i,j}$ for 
$1 \le i < j \le n$ if we define  $y_1 = x_1$. 
To show that $\theta$ determines a well-defined homomorphism we first define the special 
element $X = x_1x_2\dots x_n$ in $A(C_n)$ and 
establish the following result.
\begin{proposition}\label{Xaction}
For any reflection $R$ in $C_n$,
\[X\theta(\{R\})X^{-1} = \theta(\{\gamma R \gamma^{-1}\}).\]
\end{proposition}
\Proof{Proof. }
Since $X = x_1Y_{1,n}$ and $x_1$ commutes with $x_3, \dots , x_n$, it follows that 
$Xx_i = x_{i+1}X$ for $2 \le i < n$ and 
$Xx_{i,j} = x_{i+1, j+1}X$ for $1 \le i < j < n$.
This establishes the proposition for $R$ of the form $\prd{i,j}$ for $1 \le i < j < n$.
\vskip .2cm
The identity $Xy_j = y_{j+1}X$ for $1 \le j <n$ is a consequence of the following calculation.
\begin{eqnarray*}
Y_{2,j+1}Xy_j &=& x_2Y_{3,j+1}Xy_j = x_2XY_{2,j}y_j = x_2Xx_1Y_{2,j} \\
&=& x_2x_1x_2Y_{3,n}x_1Y_{2,j} = x_2x_1x_2x_1Y_{3,n}Y_{2,j}\\
&=& x_1x_2x_1x_2Y_{3,n}Y_{2,j} = x_1x_2XY_{2,j} = x_1x_2Y_{3,j+1}X \\
&=& x_1Y_{2,j+1}X = Y_{2,j+1}y_{j+1}X
\end{eqnarray*}
This establishes the proposition for $R$ of the form $[j]$ for $1 \le i < n$.
\vskip .2cm
Conjugating $y_n$ by $X$ gives $x_1$, since
\[Xy_n = (x_1x_2\dots x_n)y_n = x_1(x_2\dots x_ny_n) = x_1(x_1\dots x_n). \]
This establishes the proposition for the reflection $[n]$.
\vskip .2cm
Next we show $Xx_{i,n} = z_{1,i+1}X$.
\begin{eqnarray*}
z_{1,i+1}X &=& z_{1,i+1}x_1Y_{1,n} = x_1x_{1,i+1}Y_{1,n} = 
x_1x_{1,i+1}Y_{1,i}Y_{i,n} \\
&=& x_1Y_{1,i}x_{i+1}Y_{i,n}
+  x_1Y_{1,i}x_{i+1}Y_{i+1,n}x_{i,n} = x_1Y_{1,n}x_{i,n} = Xx_{i,n}
\end{eqnarray*}
This establishes the proposition for $R$ of the form $\prd{i,n}$ for $1 \le i < n$.
\vskip .2cm
The identity $Xz_{i,j} = z_{i+1,j+1}X$ for $1 \le i < j < n$ follows from the 
definition of $z_{i,j}$ and the corresponding identities for $x_{i,j}$ and $y_i$, 
which establishes the proposition for $R$ of the form $\prd{i,-j}$ for $1 \le i<j < n$.
\vskip .2cm
Next we observe that, for $3 \le j \le n$, $y_jz_{1,j} = x_{1,j}y_j$ because 
\begin{eqnarray*}
Y_{1,j}y_jz_{1,j}x_1 &=& x_1Y_{1,j}z_{1,j}x_1 = x_1Y_{1,j}x_1x_{1,j} 
= x_1x_2Y_{2,j}x_1x_{1,j} \\
&=& x_1x_2x_1Y_{2,j}x_{1,j} = x_1x_2x_1Y_{1,j} =  x_1x_2x_1x_2Y_{2,j}\\
 &=& x_2x_1x_2x_1Y_{2,j} = x_2x_1x_2Y_{2,j}x_1 = x_2x_1Y_{1,j}x_1 \\
 &=& x_2Y_{1,j}y_jx_1 = Y_{1,j}x_{1,j}y_jx_1.
\end{eqnarray*}
Since $Xz_{i,n}y_i = Xy_ix_{i,n} = y_{i+1}z_{1,i+1}X = x_{1,i+1}y_{i+1}X = x_{1,i+1}Xy_{i}$, 
it follows that $Xz_{i,n} = x_{1,i+1}X$ and hence the proposition is established 
for the final case, $R$ of the form $\prd{i,-n}$ for $1 \le i < n$.  
\hfill \Eop
\begin{defn} We define a lift of $\gamma$ to $A(C_n)$ to be an element of the form 
\[E = \theta(\{R_1\})\theta(\{R_2\})\dots \theta(\{R_n\}), \]
where the $R_i$ are reflections in $C_n$ satisfying $R_1R_2\dots R_n = [1,2,3,\dots , n]$. 
\end{defn}
\vskip .2cm
We note that one lift of $\gamma$ to $A(C_n)$ is 
\[X  = x_1x_2\dots x_n = \theta(\{[1]\})\theta(\{\prd{1,2}\})\dots \theta(\{\prd{n-1,n}\}).\]
To show that $\theta$ is well-defined it suffices, by 
Lemma~\ref{newpres}, to prove the following.
\begin{prop}
For any lift $E$ of $\gamma$ to $A(C_n)$ we have $E = X$.
\end{prop}
\Proof{Proof.} Given a lift $E = \theta(\{R_1\})\theta(\{R_2\})\dots \theta(\{R_n\})$ 
of $\gamma$ to $A(C_n)$, 
we know that $R_1R_2 \dots R_n = [1,2,\dots , n]$ and by Lemma~\ref{balancedlength} 
exactly one of the $R_k$ is of the form $[j]$.  Since $E = X$ if and only if 
$X^lEX^{-l} = X$ for any integer $l$, we may assume by the previous proposition 
that $R_k = [1]$.  We will construct a new lift $E'$ of $\gamma$ satisfying $E' = E$ and
\[E' = \theta(\{R_1\})\dots \theta(\{R_{k-2}\})\theta(\{[1]\})\theta(\{R'\})
\theta(\{R_{k+1}\})\dots \theta(\{R_n\}),\]
for some reflection $R'$.
\vskip .2cm
To simplify notation we set $R_{k-1} = T$ so that $R_{k-1}R_k = T[1]$.  
Since $T[1] \le \gamma$ we know that $T \le \gamma [1]$ or 
\[T \le \prd{1, - 2, - 3, \dots , - n}\] 
so that $T$ has the form $\prd{1, - p}$ for $2 \le p \le n$ or 
$T$ has the form $\prd{i,j}$ with $2 \le i <j \le n$.  In the latter case $\theta(\{T\})$ lies 
in the subgroup of $A(C_n)$ generated by $\{x_3, x_4, \dots , x_n\}$ and so commutes with 
$\theta(\{[1]\}) = x_1$.  Thus we can use $R' = T$.   In the former case, 
$\theta(\{T\}) = z_{1,p}$ and $E'$ can be constructed using 
\[\theta(\{T\})\theta(\{[1]\}) = z_{1,p}x_1 = x_1x_{1,p} = 
\theta(\{[1]\})\theta(\{\prd{1,p}\}).\]
After $k-1$ such steps we get 
$E = x_1\theta(\{S_2\})\dots \theta(\{S_n\})$, where the product on the right is 
a lift of $\gamma$ to $A(C_n)$. 
However, this means $S_2S_3\dots S_n = \prd{1,2,\dots ,n}$ in $C_n$ so that 
$S_i \in \Sigma_n < C_n$ and \[\theta(\{S_2\})\dots \theta(\{S_n\}) = x_2x_3 \dots x_n,\]
by Lemma~4.6 of \cite{brady1}. \hfill \Eop
\vskip .2cm
Combining the results in this subsection we get the following theorem.
\begin{thm}  The poset group $\Gamma(C_n, \gamma)$ is isomorphic to the Artin group 
$A(C_n)$ for $\gamma$ a Coxeter element in $C_n$.
\end{thm}

\subsection{The $D_n$ case}  In this case our approach will be exactly as in the $C_n$ case.  
However, the computations are more numerous and more complicated.  
The Artin group $A(D_n)$ has a presentation with generating set 
$x_1, x_2, \dots x_n$, 
subject to the relations 
\begin{eqnarray*}
x_1x_2 &=& x_2x_1, \\
x_1x_3x_1 &=& x_3x_1x_3,\\ 
x_1x_i &=& x_ix_1, \quad \mbox{ for }\quad i \ge 4 \\
x_ix_{i+1}x_i &=& x_{i+1}x_ix_{i+1}, \quad \mbox{ for }\quad 1 < i < n\quad \mbox{ and }\\
x_ix_j &=& x_jx_i, \quad \mbox{ for }\quad |j-i| \ge 2\quad \mbox{ and }\quad i,j \ne 1.
\end{eqnarray*}

\begin{defn}  We define a function $\phi$ from the generators of $A(D_n)$ to $\Gamma(D_n, \gamma)$ by 
\[x_1 \mapsto \{\prd{1, - 2}\}, x_2 \mapsto \{\prd{1,2}\}, x_3 \mapsto \{\prd{2,3}\}, \dots , 
x_n \mapsto \{\prd{n-1, n}\}\]
\end{defn}
\begin{lem}The function $\phi$ determines a well-defined surjective homomorphism.
\end{lem}
\emph{Proof:}  \hspace{.2cm} The relations involving $\phi(x_1)$ hold in $\Gamma(D_n,\gamma)$ by virtue 
of the following identities in $\Gamma(D_n, \gamma)$. 
\[\{\prd{1, - 2}\}\{\prd{1,2}\} =  \{[1][2]\} = \{\prd{1,2}\}\{\prd{1, - 2}\}\]
\begin{eqnarray*}
\{\prd{1, - 2}\}\{\prd{2,3}\}\{\prd{1, - 2}\} &=& \{\prd{1, - 2, - 3}\}\{\prd{1, - 2}\}\\
&=& \{\prd{2,3}\}\{\prd{1, - 3}\}\{\prd{1, - 2}\}\\
&=& \{\prd{2,3}\}\{\prd{1, - 2, - 3}\}\\
&=& \{\prd{2,3}\}\{\prd{1, - 2}\}\{\prd{2,3}\}
\end{eqnarray*}
\[\{\prd{1, - 2}\}\{\prd{i, i+1}\} = \{\prd{i, i+1}\}\{\prd{1, - 2}\}, 
\quad \mbox{for} \quad i\ge 3. \]
The image of the subgroup generated by $\{x_2, \dots , x_n\}$ again lies in the copy of the 
braid group corresponding to $\Sigma_n < D_n$ so that the relations not involving $\phi(x_1)$ 
hold by Lemma~4.2 and Lemma~4.4 of \cite{brady1}.  Thus $\phi$ is well-defined.

To establish surjectivity, note that both $\{\prd{i,j}\}$ and 
$\{\prd{i,i+1,\dots,j}\}$ lie in $\textrm{im}(\phi)$, for $1 \le i < j \le n$ 
as in the $C_n$ case.   To find the other reflection lifts in $\textrm{im}(\phi)$ first note that 
\[\phi(x_1x_2\dots x_j) = \{[1][2,3\dots ,j]\} = \{\prd{1, -2}\}\{\prd{1,2, \dots , j}\} 
\in \textrm{im}(\phi),\]
and $\{\prd{1,-j}\}\in \textrm{im}(\phi)$ for $j \ge 3$ since
\[\{\prd{1,2,\dots , j}\}\{\prd{1,- j}\} = \{[1][2,\dots, j]\}.\]
Reflection lifts of the form $\{\prd{2,-j}\}$ for $j \ge 3$ lie in $\textrm{im}(\phi)$ since
\[\{\prd{1, -2}\}\{\prd{1,j}\} = \{\prd{1,j,-2}\} = \{\prd{2, -j}\}\{\prd{1, -2}\}\]
and reflection lifts of the form $\{\prd{i,-j}\}$ for 
$3 \le i < j \le n$ lie in $\textrm{im}(\phi)$ since
\[\{\prd{i, -j}\}\{\prd{1,i}\}\{\prd{1,-i}\} = \{[1][i,j]\} = 
\{\prd{1,i}\}\{\prd{1,-i}\}\{\prd{i,j}\}.\]
\ \hfill \Eop
\vskip .2cm
To construct an inverse to $\phi$ we will use the presentation for $\Gamma(D_n, \gamma)$ 
given by Lemma~\ref{newpres}.
\begin{definition}
We define a function $\theta$ from the generators of $\Gamma(D_n, \gamma)$ to $A(D_n)$ by 
\[\{\prd{1,-2}\} \mapsto x_1,\quad  \{\prd{i,j}\} \mapsto x_{i,j} \quad \mbox{ and }
\quad \{\prd{i,-j}\} \mapsto z_{i,j},\] 
for $1 \le i < j \le n$, 
where $z_{i,j}$ is the unique element of $A(D_n)$ satisfying
\begin{eqnarray*}
z_{1,j}x_1 &=& x_1x_{2,j} \quad \mbox{ when } \quad j \ge 3\\
z_{2,j}x_1 &=& x_1x_{1,j}  \quad \mbox{ when } \quad  j \ge 3\\
z_{i,j}x_{1,i}z_{1,i} &=&  x_{1,i}z_{1,i}x_{i,j}  \quad \mbox{ when } \quad 3 \le i < j \le n
\end{eqnarray*}
\end{definition}
\vskip .2cm
We note that $z_{1,2} = x_1$.  Since each $x_{i,j}$ lies in the copy of $B_n$ generated by 
$\{x_2, \dots x_n\}$ the elements $x_{i,j}$  
satisfy the same identities as in the $C_n$ case. The homomorphism 
determined by $\theta$ will be surjective since 
each $x_i$ is the image of some reflection lift. 
To show that $\theta$ determines a well-defined homomorphism we define the special element 
$X = x_1x_2\dots x_n$ in $A(D_n)$ and 
establish the $D_n$ analogue of Proposition~\ref{Xaction} .
\begin{proposition}  For any reflection $R$ in $D_n$,
\[X\theta(\{R\})X^{-1} = \theta(\{\gamma R \gamma^{-1}\}).\]
\end{proposition}
\Proof{Proof. }
Since $X = x_1Y_{1,n}$ and $x_1$ commutes with $x_4, \dots , x_n$ it follows that 
$Xx_i = x_{i+1}X$ for $3 \le i < n$ and 
$Xx_{i,j} = x_{i+1, j+1}X$ for $3 \le i < j < n$.  This establishes the proposition 
in the case $R = \prd{i,j}$ for  $3 \le i < j < n$.
\vskip .2cm
For some of the later cases we will require the identities $x_{2,j}z_{1,j} = x_1x_{2,j}$ 
and $x_{1,j}z_{2,j} = x_1x_{1,j}$ for $3\le j \le n$.  The first follows from
\begin{eqnarray*}
Y_{3,j}x_{2,j}z_{1,j}x_1 &=& x_3Y_{3,j}z_{1,j}x_1 
= x_3Y_{3,j}x_1x_{2,j} = x_3x_1Y_{3,j}x_{2,j}\\
&=& x_3x_1x_{3}Y_{3,j} = x_1x_{3}x_1Y_{3,j} = 
x_1x_{3}Y_{3,j}x_1 = x_1Y_{3,j}x_{2,j}x_1 \\
&=&  Y_{3,j}x_1x_{2,j}x_1 ,
\end{eqnarray*}
while the second follows from 
\begin{eqnarray*}
x_1Y_{2,j}x_{1,j}z_{2,j}x_1 &=& x_1x_2Y_{2,j}z_{2,j}x_1 = x_1x_2Y_{2,j}x_1x_{1,j}
= x_1x_2x_3Y_{3,j}x_1x_{1,j}\\
&=& x_1x_2x_3x_1Y_{3,j}x_{1,j} = x_2x_1x_3x_1Y_{3,j}x_{1,j} = x_2x_3x_1x_3Y_{3,j}x_{1,j}\\
&=& x_2x_3x_1Y_{2,j}x_{1,j} = x_2x_3x_1x_2Y_{2,j} = x_2x_3x_1x_2Y_{2,j}\\
&=& x_2x_3x_2x_1Y_{2,j} = x_3x_2x_3x_1Y_{2,j} = x_3x_2x_3x_1x_3Y_{3,j}\\
&=& x_3x_2x_1x_3x_1Y_{3,j} = x_3x_2x_1x_3Y_{3,j}x_1 = x_3x_2x_1Y_{2,j}x_1 \\
&=& x_3x_1x_2Y_{2,j}x_1 = x_3x_1Y_{2,j}x_{1,j}x_1 = x_3x_1x_3Y_{3,j}x_{1,j}x_1\\
&=& x_1x_3x_1Y_{3,j}x_{1,j}x_1 = x_1x_3Y_{3,j}x_1x_{1,j}x_1.
\end{eqnarray*}
The conjugation action of $X$ on $x_1$ is given by $Xx_1 = x_{1,3}X$ since 
\begin{eqnarray*}
x_3Xx_1 &=& x_3x_1x_2x_3Y_{3,n}x_1 = x_3x_1x_2x_3x_1Y_{3,n} = x_3x_2x_1x_3x_1Y_{3,n}\\
&=& x_3x_2x_3x_1x_3Y_{3,n} = x_2x_3x_2x_1x_3Y_{3,n} = x_2x_3x_1x_2x_3Y_{3,n} \\
&=& Y_{1,3}X = x_3x_{1,3}X.
\end{eqnarray*}
A similar calculation gives $x_3Xx_2 = x_1x_3X$.  Since 
\[x_1x_3X = x_1x_{2,3}X = x_{2,3}z_{1,3}X\] 
we get $Xx_2 = z_{1,3}X$.  
This establishes the proposition 
in the cases $R = \prd{1,-2}$ and $R =  \prd{1,2}$.
\vskip .2cm
Next we establish $Xx_n = z_{2,n}X$. 
\[Xx_n = x_1Y_{1,n}x_n = x_1x_{1,n}Y_{1,n-1}x_n = z_{2,n}x_1Y_{1,n} = z_{2, n}X\]
which takes care of the case $R = \prd{n-1,n}$.
To obtain the identity $Xx_{1,j} = z_{1,j+1}X$ we note that
\[Y_{1,n}x_{1,j}Y_{1,j-1} = Y_{1,n}Y_{1,j} = Y_{2,j+1}Y_{1,n} = x_{2,j+1}Y_{2,j}Y_{1,n} = 
x_{2,j+1}Y_{1,n}Y_{1,j-1}\]
giving $Y_{1,n}x_{1,j} = x_{2,j+1}Y_{1,n}$ so that
\[Xx_{1,j} = x_1Y_{1,n}x_{1,j} = x_1x_{2,j+1}Y_{1,n} = z_{1,j+1}x_1Y_{1,n} = z_{1,j+1}X.\]
This completes the case $R = \prd{1,j}$ for $2 \le j < n$.
\vskip .2cm
For the identity $Xx_{1,n} = x_2X$ we compute
\[Xx_{1,n} = x_1x_2(x_3\dots x_n)x_{1,n} = x_1x_2(x_2x_3\dots x_n) = x_2X,\]
which establishes the case $R = \prd{1,n}$.  
\vskip .2cm
For $2 \le i < n$ we have 
\begin{eqnarray*}
Xx_{i,n} &=& x_1Y_{1,i+1}Y_{i+1,n}x_{i,n} = x_1Y_{1,i+1}x_{i+1}Y_{i+1,n}\\
&=& x_1x_{1,i+1}Y_{1,i}x_{i+1}Y_{i+1,n} = z_{2,i+1}x_1Y_{1,n} = z_{2,i+1}X
\end{eqnarray*}
and hence the proposition is true for $R = \prd{i,n}$ with $2 \le i < n$.
\vskip .2cm
The identity  $Xz_{1,j} = x_{1,j+1}X$ for $3 \le j < n$ follows from
\[Xz_{1,j}x_1 = Xx_1x_{2,j} = x_{1,3}x_{3,j+1}X = x_{1,j+1}x_{1,3}X = x_{1,j+1}Xx_{1},\]
while the identity $Xz_{1,n} = z_{1,2}X = x_1X $ follows from  
\[Xz_{1,n}x_1 = Xx_1x_{2,n} = x_{1,3}z_{2,3}X = x_1x_{1,3}X = x_1Xx_1.\]
This establishes the proposition for $R = \prd{1, -j}$ with $2 \le j \le n$.
\vskip .2cm
The identity $Xz_{i,n} = x_{2,i+1}X $ for $2\le i < n$ follows from  
\begin{eqnarray*}
Xz_{i,n}x_{1,i}z_{1,i} &=& Xx_{1,i}z_{1,i}x_{i,n} = z_{1,i+1}x_{1,i+1}z_{2,i+1}X\\
&=& z_{1,i+1}x_{1}x_{1,i+1}X = x_{1}x_{2,i+1}x_{1,i+1}X\\
&=& x_{2,i+1}z_{1,i+1}x_{1,i+1}X = x_{2,i+1}Xx_{1,i}z_{1,i}.
\end{eqnarray*}
This establishes the proposition for $R = \prd{i, -n}$ with $2 \le i < n$.
\vskip .2cm
Finally we note that  $x_{1,i}z_{1,i} = z_{1,i}x_{1,i}$ since 
\begin{eqnarray*}
x_{2,i}x_{1,i}z_{1,i} &=& x_2x_{2,i}z_{1,i} = x_2x_1x_{2,i} = x_1x_2x_{2,i}\\ 
&=& x_1x_{2,i}x_{1,i} = x_{2,i}z_{1,i}x_{1,i}.
\end{eqnarray*}
From this we deduce that $Xz_{i,j} = z_{i+1, j+1}X$ for $2 \le i <j <n$ since 
\begin{eqnarray*}
Xz_{i,j}x_{1,i}z_{1,i} &=& Xx_{1,i}z_{1,i}x_{i,j} = z_{1,i+1}x_{1,i+1}x_{i+1,j+1}X\\
 &=& z_{i+1,j+1}z_{1,i+1}x_{1,i+1}X = z_{i+1,j+1}Xx_{1,i}z_{1,i}.
 \end{eqnarray*}
This establishes the proposition for the remaining cases $R = \prd{i, -j}$ with 
$2 \le i < j < n$.
\ \hfill \Eop
\begin{defn} We define a lift of $\gamma$ to $A(D_n)$ to be an element of the form 
\[E = \theta(\{R_1\})\theta(\{R_2\})\dots \theta(\{R_n\}), \]
where the $R_i$ are reflections in $D_n$ satisfying $R_1R_2\dots R_n = [1][2,3,\dots , n]$. 
\end{defn}
\vskip .2cm
We note that one lift of $\gamma$ to $A(D_n)$ is 
\[X  = x_1x_2\dots x_n = \theta(\{\prd{1,-2}\})
\theta(\{\prd{1,2}\})\dots \theta(\{\prd{n-1,n}\}).\]
To show that $\theta$ determines a well-defined homomorphism it suffices, by 
Lemma~\ref{newpres}, to prove the following.
\begin{proposition}
For any lift $E$ of $\gamma$ to $A(D_n)$ we have $E = X$.
\end{proposition}
\emph{Proof:}  \hspace{.2cm}
Given a lift $E$ of $\gamma$ to $A(D_n)$, where 
\[E = \theta(\{R_1\})\theta(\{R_2\})\dots \theta(\{R_n\}),\] 
we know that $R_1R_2 \dots R_n = [1][2,\dots , n]$.  It follows for the 
proof of proposition~\ref{b_i} that one of the $R_k$ 
is of the form $\prd{1,\pm j}$.  Since $E = X$ if and only if $X^lEX^{-l} = X$ for any integer $l$, 
we may assume $R_k = \prd{1, \pm  2}$.  We treat these two cases separately.
\vskip .2cm
Suppose that $R_k = \prd{1, - 2}$.  We will construct a new lift $E'$ of $\gamma$ 
satisfying $E' = E$ and 
\[E' = \theta(\{R_1\})\dots \theta(\{R_{k-2}\})\theta(\{\prd{1, -2}\})\theta(\{R'\})
\theta(\{R_{k+1}\})\dots \theta(\{R_n\}),\]
for some reflection $R'$.
\vskip .2cm
To simplify notation we set $R_{k-1} = T$ so that 
$R_{k-1}R_k = T\prd{1,-2}$.
Since $T\prd{1, -2} \le [1][2,\dots , n]$ we know that 
\[ T \le \prd{1, -3, -4, \dots , -n, 2}\] 
so that $T$ has one of the forms 
\begin{enumerate}
\item $\prd{1,2}$,
\item $\prd{i,j}$ for $3 \le i < j \le n$, 
\item $\prd{1,-p}$ for $3 \le p \le n$ or
\item $\prd{2,-p}$ for $3 \le p \le n$. 
\end{enumerate}
\vskip .2cm
In the first case $\theta(\{T\}) = x_2$, which commutes with $\theta(\{\prd{1, -2}\}) = x_1$.  
In the second case, $\theta(\{T\}) = x_{i,j}$ lies in the subgroup generated 
by $\{x_4,\dots x_n \}$ and hence also commutes with $x_1$.   In the third case $E'$  can be 
constructed using 
\[\theta(\{T\})\theta(\{\prd{1, -2}\}) = z_{1,p}x_1 = x_1x_{2,p} = 
\theta(\{\prd{1, -2}\})\theta(\{\prd{2, p}\}) \]
and in the fourth case using 
\[ \theta(\{T\})\theta(\{\prd{1, -2}\}) = z_{2,p}x_1 = x_1x_{1,p} = 
\theta(\{\prd{1, -2}\})\theta(\{\prd{1, p}\}). \]
After $k-1$ such steps we get 
$E = x_1\theta(\{S_2\})\dots \theta(\{S_n\})$, where the product on the right is 
a lift of $\gamma$ to $A(D_n)$. However, this means 
$S_2S_3\dots S_n = \prd{1,2,\dots ,n}$ in $C_n$ so that $S_i \in \Sigma_n < C_n$ and 
\[\theta(\{S_2\})\dots \theta(\{S_n\}) = x_2x_3 \dots x_n,\]
by Lemma~4.6 of \cite{brady1}.
\vskip .2cm
Next suppose $R_k = \prd{1,  2}$. As in the previous case, we will construct a new lift 
$E'$ of $\gamma$ satisfying $E' = E$ and 
\[E' = \theta(\{R_1\})\dots \theta(\{R_{k-2}\})\theta(\{\prd{1, 2}\})\theta(\{R'\})
\theta(\{R_{k+1}\})\dots \theta(\{R_n\}),\]
for some reflection $R'$.
To simplify notation we again set $R_{k-1} = T$ so that 
$R_{k-1}R_k = T\prd{1,2}$.
Since $T\prd{1, 2} \le [1][2,\dots , n]$ we know that 
\[ T \le \prd{1, 3, 4, \dots , n, -2}\] 
so that $T$ has one of the forms 
\begin{enumerate}
\item $\prd{1,-2}$,
\item $\prd{i,j}$ for $3 \le i < j \le n$, 
\item $\prd{1,p}$ for $3 \le p \le n$ or
\item $\prd{2,-p}$ for $3 \le p \le n$. 
\end{enumerate}
\vskip .2cm
In the first case $\theta(\{T\}) = x_1$, which commutes with $\theta(\{\prd{1, 2}\}) = x_2$.  
In the second case, $\theta(\{T\}) = x_{i,j}$ lies in the subgroup generated 
by $\{x_4,\dots x_n \}$ and hence also commutes with $x_2$.   In the third case $E'$ can be 
constructed using 
\[\theta(\{T\})\theta(\{\prd{1, 2}\}) = x_{1,p}x_{1,2} = x_{1,2}x_{2,p} = 
\theta(\{\prd{1, 2}\})\theta(\{\prd{2, p}\}).\] 
In the fourth case $E'$ is constructed using 
\[\theta(\{T\})\theta(\{\prd{1, 2}\}) = z_{2,p}x_{2} = x_2z_{1,p} = 
\theta(\{\prd{1, 2}\})\theta(\{\prd{1, -p}\}).\]
The middle equality holds since
\[z_{2,p}x_2x_1 = z_{2,p}x_1x_2 = x_1x_{1,p}x_2 = 
x_1x_2x_{2,p} = x_2x_1x_{2,p} = x_2z_{1,p}x_1.\] 
\vskip .2cm
After $k-1$ such steps we   
get $E = x_2\theta(\{S_2\})\dots \theta(\{S_n\})$, where the product on the right is 
a lift of $\gamma$ to $A(D_n)$.  However, this means 
$S_2S_3\dots S_n = \prd{1,-2,\dots ,-n}$ in $C_n$ so that $S_i$ lie in the copy of $\Sigma_n$ 
generated $\{\prd{1,-2},\prd{2,3}, \dots , \prd{n-1,n} \}$ and 
\[\theta(\{S_2\})\dots \theta(\{S_n\}) = x_1x_3 \dots x_n,\]
by Lemma~4.6 of \cite{brady1}.  Finally 
\[E = x_2x_1x_3\dots x_n = x_1x_2x_3\dots x_n.  \]
\ \hfill \Eop
\vskip .2cm
Combining the results in this subsection we get the following theorem.
\begin{thm}  The poset group $\Gamma(D_n, \gamma)$ is isomorphic to the Artin group 
$A(D_n)$ for $\gamma$ a Coxeter element in $D_n$.
\end{thm}

\end{document}